\documentclass[11pt]{scrartcl}
\pdfoutput=1
\usepackage{mathptmx,graphicx,amsmath,amssymb}
\usepackage{url,color}
\usepackage[pdftex,colorlinks]{hyperref}
\definecolor{darkblue}{cmyk}{1,0,0,0.8}
\definecolor{darkred}{cmyk}{0,1,0,0.7}
\hypersetup{anchorcolor=black, citecolor=darkblue, filecolor=darkblue, menucolor=darkblue,pagecolor=darkblue,urlcolor=darkblue,linkcolor=darkblue}
\usepackage{microtype}
\usepackage{typearea}

\newcommand{\be}{\begin{equation}}
\newcommand{\ee}{\end{equation}}

\newcommand{\bea}{\begin{eqnarray}}
\newcommand{\eea}{\end{eqnarray}}
\newcommand{\R}{\mathbb{R}}
\newcommand{\Z}{\mathbb{Z}}

\newcommand{\rg}{\operatorname{rg}}
\newcommand{\id}{I}

\renewcommand{\epsilon}{\varepsilon}
\renewcommand{\phi}{\varphi}

\newtheorem{Theorem}{Theorem}
\newtheorem{Lemma}[Theorem]{Lemma}

\newtheorem{Definition}[Theorem]{Definition}

\newtheorem{Corollary}[Theorem]{Corollary}
\title{Dynamics of symmetric dynamical systems with delayed switching} 

\author{ J.~Sieber\footnote{\textit{Corresponding Author.} Centre for
    Applied Dynamics Research, School of Engineering, Fraser Noble
    Building, King's College, University of Aberdeen, Aberdeen AB24
    3UE, U.K.},\quad  
P.~Kowalczyk\footnote{ Mathematics Research Institute, Harrison Building,
  University of Exeter, Exeter, EX4~4QF, U.K.},\quad
S.J.~Hogan\footnote{Department of Engineering Mathematics, 
  University of Bristol, BS8 1TR, U.K.},\quad
M.~di~Bernardo\footnotemark[3] 
}
\begin{document}
\maketitle
\begin{abstract}
  \noindent We study dynamical systems that switch between two
  different vector fields depending on a discrete variable and with a
  delay.  When the delay reaches a problem-dependent critical value
  so-called event collisions occur.  This paper classifies and
  analyzes event collisions, a special type of discontinuity induced
  bifurcations, for periodic orbits. Our focus is on event collisions
  of symmetric periodic orbits in systems with full reflection
  symmetry, a symmetry that is prevalent in applications.  We derive
  an implicit expression for the Poincar{\'e} map near the colliding
  periodic orbit. The Poincar{\'e} map is piecewise smooth,
  finite-dimensional, and changes the dimension of its image at the
  collision. In the second part of the paper we apply this general
  result to the class of unstable linear single-degree-of-freedom
  oscillators where we detect and continue numerically collisions of
  invariant tori.  Moreover, we observe that attracting closed
  invariant polygons emerge at the torus collision.
\end{abstract}

\noindent\textbf{Keywords}: delay, relay, hysteresis, invariant torus
collision

\section{Introduction}

Relay control systems, which can be regarded as hybrid systems, are
applied in many different areas of engineering applications. Relay
feedback control might involve control of stationary processes in
industry as well as control of moving objects, for instance in flight
control. Hence, it is not surprising that in recent years much
research effort has been spent on investigations of the dynamics of
relay systems characterized by idealized off/on relay control, see for
instance \cite{DJV01,JRA99} among other works. However, relay systems
often feature intrinsic hysteretic behavior \cite{MPB03} as well as
delay in the control input. For instance, if we control a plant using
simple on/off relay feedback control via a computer network, our
control input may be delayed, due to the traffic along the network.

In the current paper we study the dynamics of relay control systems
with hysteretic behavior and delay in the control input. We wish to
understand if and how the interplay between the discontinuous
nonlinearities (switching events), hysteresis and delay in the
switching function can lead to new bifurcations. We expect so-called
\emph{event-collisions} to occur as observed in \cite{CDHK06}. We show
that they cause a change in the dimension of the phase space that
describes the local dynamics of the system, which leads to interesting
dynamical phenomena such as closed invariant polygons and collisions of
smooth invariant tori. 

\begin{figure}[t]\centering
\includegraphics[scale=0.55]{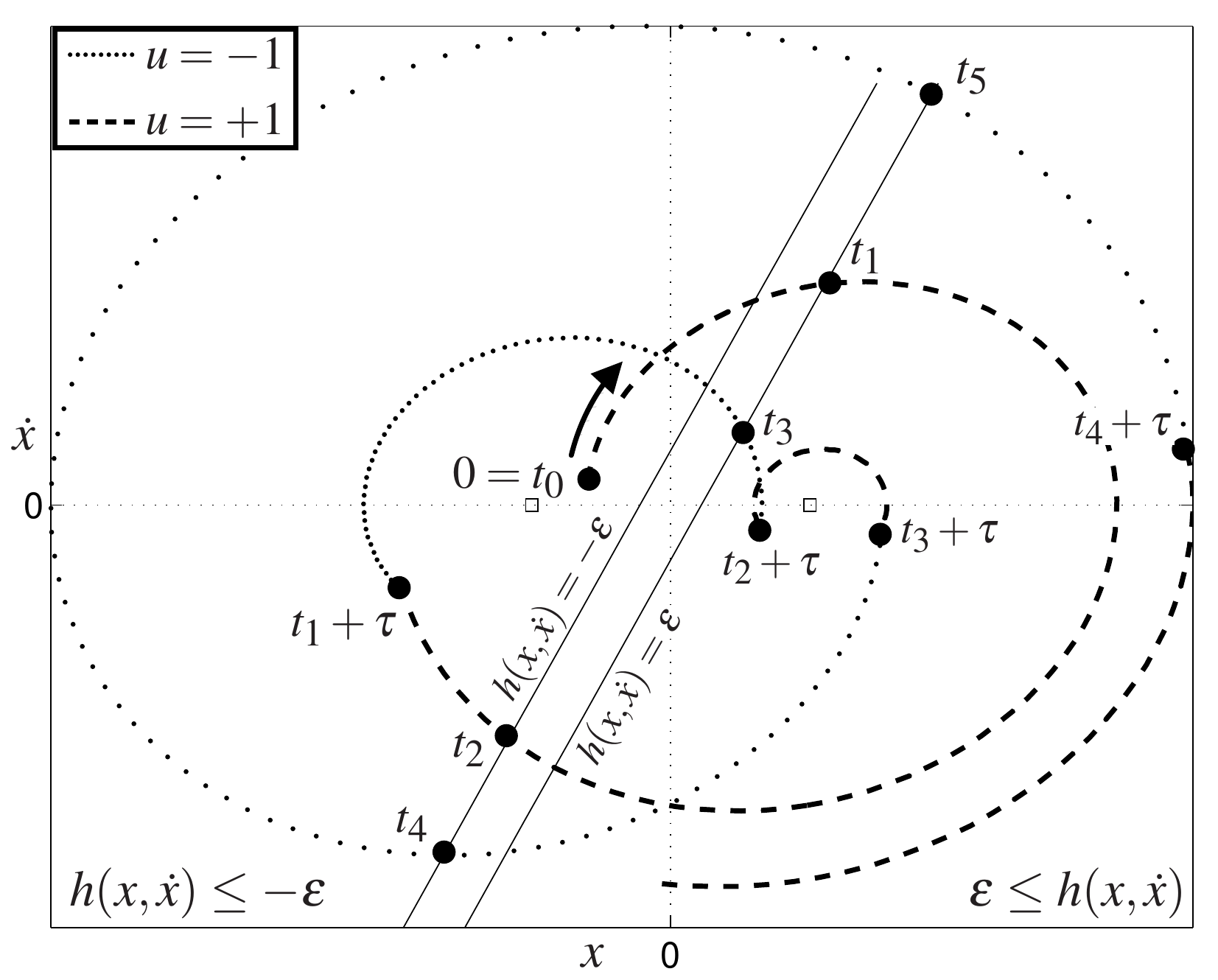}
\caption{Sketch of qualitative behavior of hybrid systems with delay.
  Initial history segment of $(x,\dot x)$ is assumed to be constant,
  initial value of $u$ is $+1$. The solid lines are the switching
  manifolds $h=\pm \epsilon$. At the times $t_j$ the trajectory
  crosses the switching line, at the times $t_j+\tau$ it switches the
  control $u$ to $-u$. }
\label{fig:sketch}
\end{figure}

Let us explain the problem for the classical prototype example of an
oscillator subject to a relay switch:
\begin{equation}
  \label{eq:osci}
  \ddot x+\zeta\dot x+\omega^2x=u\mbox{.}
\end{equation}
This example will be studied in detail in Section~\ref{sec:osci} and
in Section~\ref{sec:tor}.  In \eqref{eq:osci} $y=[x,\dot x]^T$ is the
continuous variable and $u$ provides discrete feedback in the form of
a relay. That is, $u$ is ruled by the switching law:
\begin{equation}\label{eq:osci:disc}
  u(t)=
  \begin{cases}
    \mp1 &\mbox{if $\pm h^Ty(t-\tau)\geq\epsilon$}\\
    \lim_{s \nearrow t} u(s) &\mbox{if 
      $h^Ty(t-\tau)\in(-\epsilon,\epsilon)$.}
  \end{cases}
\end{equation}
Roughly speaking, the switching law \eqref{eq:osci:disc} means that we
set $u$ to $-1$ whenever we observe that $h^Ty\geq\epsilon$ and we set
$u$ to $+1$ whenever we observe that $h^Ty\leq-\epsilon$ (negative
feedback).  The region $\{|h^Ty|<\epsilon\}$ provides a hysteresis in
the relay (a buffer between subsequent switchings from $u$ to $-u$):
we leave $u$ at its value when we observe that $|h^Ty|<\epsilon$. The
lines $\{y:h^Ty=\pm\epsilon\}$ are called the \emph{switching lines}
or, more generally for higher dimensional $y$ and nonlinear switching
functions $h$, \emph{switching manifolds}. An ordinary differential
equation such as \eqref{eq:osci} with a switching law of the form
\eqref{eq:osci:disc} is one of the simplest forms of \emph{hybrid
  dynamical systems} \cite{LJSZS03}. The twist in problem
\eqref{eq:osci},\,\eqref{eq:osci:disc} comes from the presence of a
delay $\tau$ in the observation (or the switching) in
\eqref{eq:osci:disc}.  Figure~\ref{fig:sketch} illustrates how a
trajectory of \eqref{eq:osci} might look like. The presence of a small
delay $\tau$ is often equivalent to a perturbation of the switching
law (only if $\epsilon$ is positive, see, for example, \cite{CDHK06}).  For
large delays $\tau$ the dynamics shows an enormous amount of
complexity \cite{BKW06,BH98,CDHK06,H91} because the number of
switching times in the interval $[t-\tau,t]$ is effectively included
into the dimension of the system which then becomes larger than the
dimension of the physical space \cite{S06}.

This paper studys the transition from `small' to `large' delay.  This
transition is characterized by event collisions, which are
\emph{discontinuity induced bifurcations} \cite{CDHK06}.  The notion
of a discontinuity induced bifurcation can be made precise for
periodic orbits of \eqref{eq:osci},\,\eqref{eq:osci:disc}.  The
possible types of codimension one discontinuity induced bifurcations
for generic systems with delayed switching have been classified and
analyzed in \cite{S06} (without hysteresis):
{\renewcommand{\theenumi}{\Alph{enumi}}
\begin{enumerate}
\item\label{int:gengrazing} \textbf{(grazing)} the periodic orbit touches the
  switching manifold quadratically without intersecting it,
\item\label{int:gencollision} \textbf{(corner/event collision)} the continuous
  variable $y$ is on the switching manifold at two times with
  difference equal to the delay $\tau$: say, at time $t$ and time
  $t-\tau$. The name is due to the fact that the periodic orbit
  typically has a `corner' at $t$ because it switches the vector
  field ($u\mapsto-u$) at $t$.
\end{enumerate}
}\noindent Most of the event collisions reported in \cite{CDHK06}
violate one of the secondary genericity conditions in \cite{S06},
namely that only one corner collides.  The reason behind this
unexpected degeneracy of the phenomena observed in \cite{CDHK06} is
that the oscillator system~\eqref{eq:osci},\,\eqref{eq:osci:disc}
(which is the subject of study in \cite{CDHK06}) is symmetric with
respect to reflection at the origin:
\begin{equation}
  y\mapsto-y\mbox{,\quad} u\mapsto-u\mbox{.}\label{eq:intsym}
\end{equation}
The periodic orbits of primary interest inherit this symmetry. If the
control switches from $u$ to $-u$ at time $t$ it switches from $-u$ to
$u$ one half-period later (at $t+p/2$ if $p$ is the period).  Whenever
the corner at $t$ is on the switching manifold the opposite corner (at
$t+p/2$) is on the switching manifold as well due to reflection
symmetry. This is in contrast to the secondary genericity conditions
in \cite{S06}. In fact, all systems studied in
\cite{BKW06,BH98,CDHK06,H91} have a reflection symmetry because they
are piecewise linear, which renders the generic theory inapplicable to
symmetric periodic orbits in these examples of practical interest.

We close this gap in the general context, that is, for $y\in\R^n$ and
nonlinear right-hand-sides and switching laws.  The first main result
of the paper is an implicit expression for the local return map for
symmetric periodic orbits close to a simultaneous collision of two
corners. The simplest example for this scenario is a symetric periodic
orbit of period $p=2\tau$ for which $y(0)$ and $y(\tau)$ lie both on
the switching manifold. Due to the reflection symmetry this event has
codimension one. The return map near a colliding orbit is a piecewise
smooth map $F$ that has a $(n-1)$-dimensional image in one subdomain
of its phase space and a $n$-dimensional image in the other subdomain.
This type of map has been studied extensively in the context of
grazing bifurcations in Filippov systems (systems where
$\epsilon=\tau=0$; see the textbook \cite{BBCK07} and the recent
review \cite{KDCHHPKN06}). In this sense our paper provides a link
between the theory of systems with delayed switching and hysteresis to
the theory of low-dimensional piecewise smooth maps as treated in \cite{BBCK07}.

The second aspect of the paper is a study of the concrete class of
oscillators \eqref{eq:osci},\,\eqref{eq:osci:disc} near event
collisions of symmetric periodic orbits. We pay special attention to
parameter regimes where the event collision leads to quasi-periodicity
and, thus, collisions of invariant tori and the appearance of closed
invariant polygons.  None of the phenomena described in the sections
\ref{sec:osci} and \ref{sec:tor} are present in \cite{CDHK06} because
the authors restricted their study to the special case of pure
position feedback $h^T=(1,0)$ in the switching law
\eqref{eq:osci:disc}. 

The first part of our paper adopts a `local' approach, which is slightly
different from the studies \cite{BKW06,BH98,CDHK06,FFS02,H91}. We do
not aim to classify the dynamics of a particular class of systems as
completely as possible.  Instead, we develop a local bifurcation
theory, considering a general system with $n$-dimensional physical
space and assuming that it has a periodic orbit $(y_*,u_*)$. Then we
study the dynamics near the periodic orbit deriving conditions for the
presence of bifurcations of a certain codimension. In this way our
results will be more general than studies of specific classes of
systems but all statements are valid only locally. The consideration
of only two vector fields and and a binary switch is not really a
restriction when one studies the local dynamics near a particular
periodic orbit. The second part of the paper then demonstrates how
this local bifurcation theory can be useful in combination with
numerical continuation to understand the dynamics of the class of
oscillators \eqref{eq:osci},\,\eqref{eq:osci:disc}.


The paper is organized as follows. 
In Section~\ref{sec:evol} we introduce the basic notation and collect
some fundamental facts about the forward evolution defined by general
systems with delayed switching and hysteresis. We also point out
differences to the neighboring cases of zero hysteresis or zero delay.
Section~\ref{sec:genper} specifies conditions on periodic orbits which
guarantee that a local return map exists and that this map is
finite-dimensional and, generically, smooth.  Section~\ref{sec:bif}
classifies and analyzes the codimension one discontinuity induced
bifurcations of periodic orbits for systems with reflection symmetry.
Section~\ref{sec:osci} studys the single-degree-of-freedom oscillator
\eqref{eq:osci},\,\eqref{eq:osci:disc} near an event collision,
classifying the local bifurcations of symmetric periodic orbits.  In
Section~\ref{sec:tor} we unfold the collision of a Neimark-Sacker
bifurcation with the switching line, a codimension two event
involving collisions of closed invariant curves and closed invariant
polygons. The appendix contains the proofs of all lemmas and the technical
details of some constructions.

\section{Fundamental properties of the evolution}
\label{sec:evol}
We consider general hybrid dynamical systems with delay of the
following form:
\begin{eqnarray}
  \dot y(t)&=&f(y(t),u(t)) \label{eq:cont}\\
  u(t)&=&
  \begin{cases}
    -1 & \mbox{if $h(y(t-\tau))\geq\epsilon$, or}\\
    &\mbox{if $h(y(t-\tau))\in(-\epsilon,\epsilon)$ and $u_-(t)=-1$,}\\
    \phantom{-}1 & \mbox{if $h(y(t-\tau))\leq-\epsilon$, or}\\
    &\mbox{if $h(y(t-\tau))\in(-\epsilon,\epsilon)$ and $u_-(t)=1$}\\
  \end{cases}\label{eq:disc}
\end{eqnarray}
where $\epsilon$ is positive.  The continuous variable $y$ is
$n$-dimensional and the discrete variable $u$ is, for simplicity of
presentation, binary, controlling the switching between the two vector
fields $Y_\pm^t$ given by $\dot y=f(y,\pm1)$.  In the definition of
$u(t)$ in (\ref{eq:disc}), $u_{-}(t)$ is defined as
\begin{equation}
  u_{-}(t):=\lim_{s \nearrow t} u(s)\mbox{.}
  \label{eq:ulimit}
\end{equation}
We assume continuous differentiability for the functions
$f(\cdot,\pm1):\R^n\mapsto\R^n$ and $h:\R^n\mapsto\R$ in the
right-hand-side of \eqref{eq:cont},\,\eqref{eq:disc} with respect to
the argument $y$ (and, possibly, parameters). Furthermore, we assume
that the gradient $h'$ is non-zero everywhere and that $f$ and $h$ have
uniform Lipschitz constants (using a prime for the
derivative with respect to $y$):
\begin{displaymath}
  \|f'\|\leq L_{\max}\mbox{, and\ }  
  \|h'\|\leq H_{\max}\mbox{.}
\end{displaymath}
Due to the delay in the argument of $y$ in the switching decision
\eqref{eq:disc} the phase space of \eqref{eq:cont},\,\eqref{eq:disc}
is infinite-dimensional.  An appropriate initial value for $y$ is the
history segment $y([-\tau,0])$ \cite{DGLW95,S89}. Thus, the phase
space of \eqref{eq:cont},\,\eqref{eq:disc} is
$C([-\Theta,0];\R^n)\times\{-1,1\}$ where $\Theta$ is an upper bound
for the delay $\tau$ (we will vary $\tau$ as a bifurcation parameter
in the sections \ref{sec:osci} and \ref{sec:tor}). The notation $C([-\Theta,0];\R^n)$ refers
to the space of continuous functions on the interval $[-\Theta,0]$
with values in $\R^n$.

We clarify in Appendix~\ref{sec:app:evolution} in which sense
system~\eqref{eq:cont},\,\eqref{eq:disc} constitutes a dynamical
system. In particular, we give a precise definition of the forward
evolution $E^t(\xi,u_0)$ from an arbitrary initial condition
$(\xi_0,u_0)\in C([-\Theta,0];\R^n)\times\{-1,1\}$ using the variation
of constants formulation of \eqref{eq:cont},\,\eqref{eq:disc}. The
evolution $E^t(\xi_0,u_0)$ has a continuous infinite-dimensional
component $E_c^t(\xi_0,u_0)\in C([-\Theta,0];\R^n)$ and a discrete
component $E_d^t(\xi_0,u_0)\in\{-1,1\}$.  The infinite-dimensional
components are related to each other by a simple time shift:
\begin{displaymath}
  E_c^t(\xi_0,u_0)(\theta)=E_c^{t+\theta}(\xi_0,u_0)(0)
\end{displaymath}
with a common \emph{headpoint trajectory} $E_c^t(\xi_0,u_0)(0)\in
C([0,t_E];\R^n)$ for all $\theta\in[-\Theta,0]$ and $t\in[0,t_E]$; see
\cite{DGLW95}.  Thus, $E_c^t(\xi_0,u_0)$ is continuous in $t$ for all
$t\geq0$. However, in general $E_c^t(\xi_0,u_0)$ (for a fixed $t$)
does not depend continuously on the component $\xi_0$ of its initial
value.

Lemma~\ref{thm:basic} summarizes two basic facts about the evolution
$E^t(\xi,u)$ that are proved in Appendix~\ref{sec:app:evolution}.
\begin{Lemma}[Fundamental properties of evolution]\label{thm:basic}
  \ \\ Let $E^t(\xi_0,u_0)$ be a trajectory of the dynamical
  system~\eqref{eq:cont},\,\eqref{eq:disc} on a bounded interval
  $[t_0,t_E]$. Then the following holds.
  \begin{enumerate}
  \item\label{thm:swfin} The discrete component $E_d^t(\xi_0,u_0)$
    changes its value (switches) only finitely many times in
    $[t_0,t_E]$. Thus, $E_c$ follows either $Y_+$ or $Y_-$ for all but
    finitely many times $t\in [t_0,t_E]$.
  \item\label{thm:swbound} If $t_0\geq\tau$ then the number $n_s$ of
    switchings of $E_d$ is bounded by
    \begin{displaymath}
      n_s\leq 1+[t_E-t_0] \cdot L_{\max}H_{\max}y_{\max}/(2\epsilon)
    \end{displaymath}
    where $y_{\max}$ is the maximum of
    $\|E_c^t(\xi_0,u_0)\vert_{t\in[t_0,t_E]}\|$ and $L_{\max}$ and
    $H_{\max}$ are the Lipschitz constants of $f$ and $h$, respectively.
  \end{enumerate}
\end{Lemma}
The statements in Lemma~\ref{thm:basic} are subtly dependent on the
presence of hysteresis and delay, which can be seen from the fact that
statement~\ref{thm:swfin} is not true in general if $\epsilon=0$.  See
\cite{S06} for details about systems with switches that have delay but
no hysteresis and \cite{BKW06,BH98} for studies of a piecewise linear
oscillator with rather intricate behavior.

Lemma~\ref{thm:basic} implies that the headpoint trajectory
$E_c^t(\xi_0,u_0)(0)$ is differentiable with respect to time $t$ in
any bounded interval $[t_0,t_E]$, following one of the flows $Y_\pm$,
except in a finite number of times.

\section{Periodic orbits}
\label{sec:genper}
If the relay is used to control a linearly unstable system, such as
the oscillator in \eqref{eq:osci} with negative damping $\zeta$, we
cannot expect to have stable equilibria. In this case the simplest
possible long-time behavior of \eqref{eq:osci} is periodic motion.  We
focus on the dynamics near periodic orbits for the remainder of the
paper. This section introduces the necessary notation and collects
some basic facts about the behavior of the general relay
system~\eqref{eq:cont},\,\eqref{eq:disc} near a periodic orbit.
\begin{Definition}[Crossing time]
  Let $L=(y_*(t),u_*(t))$ be a periodic solution of the general
  relay system \eqref{eq:cont},\,\eqref{eq:disc}, that is,
  $y_*(t+p)=y_*(t)$, $u_*(t+p)=u_*(t)$ for all $t\geq0$ and some
  period $p>0$. We call a time $t$ crossing time of $L$ if
  $|h(y_*(t))|=\epsilon$ and $u_{*,-}(t+\tau)=\lim_{s\nearrow
    t}u_*(s+\tau)\neq u_*(t+\tau)$.
\end{Definition}
Equivalently, we could say that $t$ is a crossing time if
$h(y_*(t))=\epsilon u_{*,-}(t+\tau)=\epsilon \lim_{s\nearrow
  t}u_*(s+\tau)$. Time $\tau$ after a crossing time $u_*$ will switch.
A periodic orbit $L$ can have only finitely many crossing times $t_k$
($k=1,\ldots,m$) per period. Thus, $L$ is differentiable, following
one of the flows $Y_\pm$, in all times except, possibly, in
$t_k+\tau+jp$ ($k=1,\ldots,m$, $j\in\Z$).  We assume that the period
$p$ is larger than the delay $\tau$ (without loss of generality
because $p$ does not have to be the minimal period).

First, we establish when the evolution is continuous with respect to
its initial value in a periodic orbit. We call the condition for
continuity \emph{weak transversality}. As the name suggests it
excludes that the periodic orbit touches the switching manifold
without crossing it (but does not require positive speed of crossing,
thus, \emph{weak} transversality). If this condition is satisfied then
it makes sense to define a local return (or \emph{Poincar{\'e}}) map
along the orbit
\begin{Definition}[Weak transversality]\label{def:wtran}
  Let $L=(y_*(t),u_*(t))$ be a periodic solution of period $p$
  of the general relay system \eqref{eq:cont},\,\eqref{eq:disc}.  We
  say that $L$ satisfies weak transversality if $|h(y_*(t))|$ is
  locally strictly monotone increasing near all crossing times
  $t_1,\ldots,t_m\in[0,p]$ of $L$.
\end{Definition}
Weak transversality permits that $\lim_{s\nearrow t} \dot
y_*(s)\neq\lim_{s\searrow t}\dot y_*(s)$ and each of the limits may be
tangential to the switching manifold $\{h=\pm\epsilon\}$ in any
crossing time $t$ of $L$.  However, it enforces that the periodic
orbit $L$ cannot just touch the switching manifold quadratically, and
that the number $m$ of crossing times is even. Weak transversality is
generically satisfied for a periodic orbit.
\begin{Lemma}[Continuity]\label{thm:continuity}
  If a periodic orbit $L=(y_*(t),u_*(t))$ is  weakly
  transversal then the continuous component $E_c^t(\xi,u)$ of the
  evolution $E$ is continuous with respect to $\xi$ in
  $(\xi,u)=(y_*(\cdot),u_*(0))$ for all $t\geq0$. 
\end{Lemma}
(See Appendix~\ref{sec:app:evolution} for proof.) Due to this
continuity it makes sense to define a local return map (also called
\emph{Poincar{\'e}} map) along a weakly transversal periodic orbit
$L=(y_*(t),u_*(t))$ to a local cross section ${\cal S}$. Assume
(without loss of generality) that neither $0$ nor $-\tau$ is a
crossing time of $L$ and that $y_*$ follows $Y_+$ near time $0$ (thus,
$u_*(0)=1$).  Then $y_*$ is differentiable in $0$. Let
$\Sigma=\{y\in\R^n:\dot y_*(0)^T[y-y_*(0)]=0\}$ be the hyperplane in
$\R^n$ orthogonal to $\dot y_*(0)$ in $y_*(0)$. We choose as
Poincar{\'e} section
\begin{equation}\label{eq:poincare}
  {\cal S}\!=\!\{(\xi,u)\in C([-\Theta,0];R^n)\times\{-1,1\}: 
  \xi(0)\in\Sigma\mbox{,\, } u=1\}\mbox{.}
\end{equation}
It is not actually necessary to choose $\Sigma$ orthogonal to $\dot
y_*(0)$. Any cross-section $\Sigma$ which is transversal to $y_*$ at
$t=0$ is admissible. Let
$-\tau<t_{*,1}<\ldots<t_{*,\mu}<0<t_{*,\mu+1}<\ldots<t_{*,m}<p-\tau$
be the crossing times of $L$ in $[-\tau,p-\tau]$. The following lemma
states that the local return map to ${\cal S}$ is in fact a map in
$\Sigma\times\R^{\mu}$ (a space of dimension $n-1+\mu$).  The notation
$U(v)$ refers to a (sufficiently small) neighborhood of a vector or
number $v$ in Lemma~\ref{thm:finmap} and throughout the paper.
\begin{Lemma}[Finite-dimensional Poincar{\'e} map]\label{thm:finmap}\ \\
  There exist a $\delta>0$ and neighborhoods $U_1(y(\cdot))\subseteq
  U_2(y(\cdot))\subset C([-\Theta,0];R^n)$ and $U(p)\subset\R$ such that the
  following holds: All initial conditions $(\xi,u)\in {\cal S}_1={\cal
    S}\cap [U_1(y(\cdot))\times\{1\}]$ have a unique return time
  $T(\xi)\in U(p)$ to ${\cal S}_2={\cal S}\cap
  [U_2(y(\cdot))\times\{1\}]$. For any initial condition in
  $(\xi,u)\in U_1(y(\cdot))\times\{1\}$ there exist unique times $t_1<\ldots<t_\mu$
  in $(-\tau,0)$ such that 
  \begin{equation}\label{eq:tdef}
    t_j=\min\{t\in[t_{*,j}-\delta,t_{*,j}+\delta]:|h(\xi(t))|=\epsilon\}
    \mbox{,\quad($j=1,\ldots,\mu$).}
  \end{equation} The local return map
  \begin{displaymath}
    P(\xi)=E^{T(\xi)}(\xi,1)
  \end{displaymath}
  depends only on $(\xi(0),t_1\ldots,t_\mu)\in\Sigma\times\R^\mu$.
\end{Lemma}
We have omitted the discrete component of the return map from $P$
because it is always $+1$.  An equivalent definition of $t_j$ is
\begin{math}
   t_j=\min\{t\in[t_{*,j}-\delta,t_{*,j}+\delta]:h(\xi(t))=\epsilon u_j\}
\end{math}
where $u_j=u_{*,-}(t_{*,j}+\tau)$.  The precise dependence of $P$ on
$(\xi(0),t_1,\ldots,t_\mu)$ is given in
Appendix~\ref{sec:app:evolution}.
\begin{Corollary}[Smooth Poincar{\'e} map for generic periodic orbits]
  \label{thm:smoothmap}
  Assume that all crossing times $t_{*,k}$ ($k=1,\ldots,m$) of the periodic
  orbit $L=(y_*,u_*)$ satisfy the following two conditions:
  \begin{enumerate}
  \item \label{cond:nocoll} \emph{(no collision)} the time
    $t_{*,k}-\tau$ is a not a crossing time of $L$, and
  \item\label{cond:trans} \emph{(smooth transversality)}
    $h'(y_*(t_{*,k}))\,\dot y_*(t_{*,k})\neq 0$.
  \end{enumerate}
  Then the Poincar{\'e} map $P$ depends smoothly on the coordinates
  $(\xi(0),t_1,\ldots,t_\mu)\in\Sigma\times\R^\mu$ of $(\xi,u)\in{\cal
    S}_1$ ($t_1,\ldots,t_\mu$ as defined by \eqref{eq:tdef}).
\end{Corollary}
The two conditions of Corollary~\ref{thm:smoothmap} guarantee that
whenever $u_*(t)$ changes its value then $y_*(t-\tau)$ follows one of
the two flows $Y_\pm$ in a neighborhood of $t-\tau$ and $\dot
y_*(t-\tau)$ points transversally through the switching manifold
$\{h=\pm\epsilon\}$. Generically, these two conditions are satisfied,
which implies that, generically, the Poincar{\'e} map of a periodic
orbit is smooth. Condition~\ref{cond:trans} is more restrictive than
weak transversality as it requires differentiability of $y_*(\cdot)$
and a non-zero time derivative of $h(y_*(\cdot))$ in $t_{*,k}$.

Corollary~\ref{thm:smoothmap} implies that the dynamics and possible
bifurcations near a periodic orbit $(y_*,u_*)$ of the general delayed
relay system are described by the theory for low-dimensional smooth
maps whenever the conditions~\ref{cond:nocoll} and \ref{cond:trans} are
satisfied. See, for example, \cite{K04} for a comprehensive textbook
on bifurcation theory for smooth systems.

\begin{Definition} [Slowly oscillating periodic orbit]\ \\
  We call a periodic orbit $L=(y_*,u_*)$ slowly oscillating if the
  distance between subsequent crossing times $t_k$ of $L$ is always greater
  than the delay $\tau$.
\end{Definition}
If we choose the Poincar{\'e} section ${\cal S}$ appropriately then
the local return map of a slowly oscillating periodic orbit $L$
satisfying the genericity conditions of Corollary~\ref{thm:smoothmap}
is simply a return map to $\Sigma$, the hyperplane of $\R^n$ in the
definition of ${\cal S}$. This has been observed in \cite{CDHK06} for
the oscillator \eqref{eq:osci},\,\eqref{eq:osci:disc}.

\section{Discontinuity induced bifurcations}
\label{sec:bif}

Let us assume that the right-hand-side $f$ and the switching function
$h$ depend on an additional parameter $\lambda$.  What happens to the
dynamics near a periodic orbit under variation of $\lambda$ (or,
alternatively, the delay $\tau$)? The previous section has established
that, as long as the conditions of Corollary~\ref{thm:smoothmap} are
satisfied, we should expect standard bifurcation scenarios such as
period doubling, saddle-node or Neimark-Sacker bifurcations (see
\cite{K04} for a classification). However, when varying the parameter
$\lambda$ we can also achieve that any of the conditions
\ref{cond:nocoll} and \ref{cond:trans} fails at special parameter
values. We call these events \emph{discontinuity induced
  bifurcations}.

For compactness of presentation we assume that for $\lambda<\lambda_0$ the
periodic orbit $L=(y_*,u_*)$ is slowly oscillating, that is, the
distance between subsequent crossing times $t_k$ is always greater
than the delay $\tau$ for $\lambda<\lambda_0$. 
\subsection{Generic bifurcations}
\label{sec:genbif}

\paragraph*{Generic grazing}
If condition~\ref{cond:trans} of Corollary~\ref{thm:smoothmap} is
violated at $\lambda=\lambda_0$  the orbit $L$ grazes the
switching manifold tangentially at a crossing time $t_*$. Let us
assume that $y_*$ follows $Y_+$ in $t_*$ without loss of generality.
That is,
\begin{displaymath}
  0=\frac{d}{dt} h(y_*(t),\lambda_0)\vert_{t=t_*}=
  h'(y_*(t_*),\lambda_0)\,\dot y_*(t_*)=
  h'(y_*(t_*),\lambda_0)\,f(y_*(t_*),+1,\lambda_0)\mbox{.}
\end{displaymath}
Generically, one can expect that
\begin{equation}\label{eq:quadgrazing}
  0\neq \frac{d^2}{(dt)^2}h(y_*(t),\lambda_0)\vert_{t=t_*}=
  h'(y_*(t_*),\lambda_0)\,\ddot y_*(t_*)+ h''(y_*(t_*),\lambda_0)\,
  [\dot y_*(t_*)]^2\mbox{,}
\end{equation}
which means that the periodic orbit $L$ touches the switching manifold
quadratically and not to a higher order. However, under condition
\eqref{eq:quadgrazing} $L$ is not weakly transversal. In general, we
cannot expect that the evolution is continuous in $L$. Thus,
trajectories arbitrarily close to $L$ leave a fixed neighborhood of
$y_*$ in a finite time (typically one period). Consequently, generic
grazing of the periodic orbit $L$ cannot be described by the approach
of local bifurcation theory adopted in this paper.

\paragraph*{Generic corner collision}
If condition~\ref{cond:trans} of Corollary~\ref{thm:smoothmap} is
violated at $\lambda=\lambda_0$ the orbit $L$ switches between
the two vector fields exactly at a crossing time $t_1$. That is,
\begin{displaymath}
u_*(t_1)\neq \lim_{s\nearrow t_1} u_*(s)=u_{*,-}(t_1)\mbox{.}
\end{displaymath}
Let us denote $f_1=f(y_*(t_1),u_{*,-}(t_1),\lambda_0)$ and
$f_2=f(y_*(t_1),u_*(t_1),\lambda_0)$. Generically, we can expect that
{\renewcommand{\theenumi}{C\arabic{enumi}}
\begin{enumerate}
\item\label{cond:nosim}  all other crossing times $t_2,\ldots,t_m$ of
  $L$ do not collide, that is, $u_*(t_k)=u_{*,-}(t_k)$ for
  $k=2,\ldots,m$, and
\item\label{cond:stran}  the
  left-sided tangent $f_1$ and the right-sided tangent $f_2$ to $y_*(\cdot)$
  in $t_1$ are transversal to the switching manifold. More precisely,
  \begin{equation}
    \label{eq:stran}
    q:=h'(y_*(t_1),\lambda_0)\,f_1\cdot 
    h'(y_*(t_1),\lambda_0)\,f_2\neq 0\mbox{.}
  \end{equation}
\end{enumerate}
} If $q<0$ in condition \eqref{eq:stran} the periodic orbit $L$ has a
corner at $t_1$ and this corner touches the switching manifold from
one side at the crossing time $t_1$. Thus, $L$ is not weakly
transversal to the switching manifold in its crossing time $t_1$.  In
this case we cannot expect that the evolution is continuous in $L$.
Consequently, corner collisions with $q<0$ cannot be described using
local bifurcation theory, either. An explicit expression for the
return map (which is a piecewise asymptotically linear
$(n-1)$-dimensional map) under the assumptions \ref{cond:nosim} and
\ref{cond:stran} for the case $q>0$ has been derived in \cite{S06} (as
case (a) in Appendix D of \cite{S06}).

\subsection{Corner collision with reflection symmetry}
\label{sec:symbif}
Often the practically relevant examples have special
symmetries which enforce that if one corner of a symmetric
periodic orbit collides at $\lambda=\lambda_0$ then other corners of the
symmetric periodic orbit collide at $\lambda_0$ simultaneously. For
example, the systems studied in \cite{BKW06,BH98,CDHK06,H91} and our
prototype oscillator \eqref{eq:osci} are all piecewise affine:
\begin{equation}
  \label{eq:lin}
  \begin{split}
    f(y,u)&=Ay+bu\\
    h(y)&=h^T y\mbox{.}
  \end{split}
\end{equation}
Thus, they have a full reflection ($\Z_2$) symmetry
\begin{equation}
  \label{eq:sym}
  \begin{split}
    f(y,u,\lambda)&=-f(-y,-u,\lambda)\\
    h(y,\lambda)&=-h(-y,\lambda)\mbox{.}
  \end{split}
\end{equation}
This means that, even though condition~\ref{cond:nosim} should be
generically satisfied for a colliding periodic orbit, restricting to the
generic case disregards many practically relevant systems.

The $\Z_2$ symmetry \eqref{eq:sym} typically gives rise to a symmetric
periodic orbit $L=(y_*,u_*)$ satisfying $y_*(t-T)=-y_*(t)$ and
$u_*(t-T)=-u_*(t)$ for the half-period $T$ and all times $t$. A corner
collision of $L$ for a crossing time $t$ at a
special parameter $\lambda=\lambda_0$ automatically induces a corner
collision for the crossing time $t-T$, a scenario that is not covered
by the generic bifurcation scenarios listed in
Section~\ref{sec:genbif}.

Let us assume that system~\eqref{eq:cont},\,\eqref{eq:disc} has full
reflection symmetry \eqref{eq:sym} and a symmetric periodic orbit
$L=(y_*,u_*)$ of half-period $T$ that experiences a corner collision
at the parameter $\lambda=\lambda_0$ for crossing time $0$ and,
enforced by symmetry, for crossing time $T$. For compactness of
presentation let us assume that $0$ and $T$ are the only crossing
times of $L$. This implies that the delay $\tau$ equals the
half-period $T$ and that $u_*$ switches between $+1$ and $-1$ at the
crossing times $0$ and $T$. Without loss of generality $L$ consists of
the two segments
\begin{alignat*}{2}
  y_*([0,T])\phantom{2}&=Y_+^{[0,T]}y_*(0)\mbox{,\quad} &
  u_*([0,T))\phantom{2}&=+1\\
  y_*([T,2T])&=Y_-^{[0,T]}y_*(T)=-Y_+^{[0,T]}y_*(0)&\mbox{,\quad} 
  u_*([T,2T))&=-1
\end{alignat*}
Moreover, $h(y_*(0),\lambda_0)=\epsilon$ and
$h(y_*(T),\lambda_0)=-\epsilon$. Note that the colliding orbit $L$ is
always following the `wrong' flow.  That is, the orbit is identical
in shape to the periodic solution with positive feedback ($+1$ and
$-1$ interchanged in \eqref{eq:disc}) and zero delay. In addition we
assume that condition~\ref{cond:stran} is satisfied for the collision
time $t=0$ (and, by symmetry, for $t=T$) with $q>0$. We call this
condition \emph{strict transversality} because it is stronger than the
weak transversality introduced in Section~\ref{sec:genper}:
\begin{equation}
  \label{eq:stran:qg0}
  q:=h_0'\,f_1\cdot h_0'f_2>0
\end{equation}
where $h_0'=h'(y_*(0),\lambda_0)$, $f_1=f(y_*(0),-1,\lambda_0)$ and
$f_2=f(y_*(0),1,\lambda_0)$. 

This strict transversality guarantees that the evolution of the
continuous component $E_c^t(\xi,u)$ is continuous for
$\xi=y_*(s+\cdot)$ for all $s\in\R$. We choose a cross section
$\Sigma$ for the Poincar{\'e} map at $t=\Delta$ ($0<\Delta\ll1$) and
orthogonal to $f_2$: $\Sigma=\{y:f_2^T[y-y_*(\Delta)]=0\}$. If
$\Delta$ is sufficiently small then $\Sigma$ is transversal to $y_*$
in $t=\Delta$ because $\dot y_*(\Delta)=f_2+O(\Delta)$.

Using the cross-section $\Sigma$ in the definition of the Poincar{\'e} map
$P$, Lemma~\ref{thm:finmap} states that the return map $P$ in its
domain of definition
\begin{displaymath} {\cal S}_1= U(y_*(\Delta+\cdot))\cap \{\xi\in
  C([-\Theta,0];\R^n):\xi(0)\in\Sigma\}
\end{displaymath}
depends for an initial value $\xi\in{\cal S}_1$ only on the headpoint
$\xi(0)\in\Sigma$ and the time
$t_1(\xi)=\min\{t\in[-2\Delta,0]:h(\xi(t))=\epsilon\}$. The time
$t_1(\xi)$ is the time when $\xi$ crosses the switching manifold
$\{h=\epsilon\}$ for the first time in $[-2\Delta,0]$. This time $t_1(\xi)$
exists if $U(y_*(\Delta+\cdot))$ is sufficiently small. Effectively,
the map $P$ depends only on $n$ coordinates.

The following lemma simplifies the representation of the Poincar{\'e}
map $P$ to a map from $U(y_*(0))$ back to $U(y_*(0))$.
\begin{Lemma}[Return map near collision]
  \label{thm:retmapreduce}
  All elements $\xi$ of the image $\rg P$ of the return map $P$ have the form
  \begin{equation}\label{eq:mapxi}
    \begin{split}
      \xi(s)&=
      \begin{cases}
        Y_+^{[s+\theta(y)]}y &\mbox{if $s\in[-\theta(y),0]$,}\\
        Y_-^{[s+\theta(y)]}y &\mbox{if $s\in[-\Theta,-\theta(y)]$,}
      \end{cases}
    \end{split}
  \end{equation}
  where $y\in U(y_*(0))$ and $\theta(y)$ is implicitly defined by the
  condition
  \begin{math}
    Y_+^{\theta(y)}y\in\Sigma\mbox{.}
  \end{math}
  The coordinates $(\xi(0),t_1(\xi))$, which are necessary for the definition of
  $P$, are uniquely defined by $y$.
\end{Lemma}
(See Appendix~\ref{sec:app:switch} for the proof.) We have omitted the
discrete variable $u$ as an argument of $P$ because it is always $+1$
at the cross-section. The point $y=\xi(-\theta(y))$ is the point where
the continuous component $\xi\in\rg P$ switches from following $Y_-$
to following $Y_+$.  The time $\theta(y)$ is the time that has elapsed
between the switching and the intersection of the headpoint with
$\Sigma$. For a sufficiently small $U(y_*(\Delta+\cdot)$ we can assume
that $\theta(y)$ is in $[0,2\Delta]$. Lemma~\ref{thm:retmapreduce}
states that the return map $P$, restricted to its image $\rg P$, can
be described as a map from the switching point of the initial value to
the switching point of its image under $P$.

Exploiting the representation \eqref{eq:mapxi} we can describe the dynamics of
the Poincar{\'e} map $P$ by a map $m$ mapping from
$U(y_*(0))\subset\R^n$ back to $U(y_*(0))$.
\begin{Theorem}[Reduced return map near collision]\label{thm:retmapform}
  Let $\lambda$ be sufficiently close to $\lambda_0$. The return map
  $m$ for elements of the image $\rg P$ of the Poincar{\'e} map $P$ is
  given by $m=F\circ F$ where $F:U(y_*(0))\mapsto U(y_*(0))$ is defined
  by
  \begin{equation}
    \label{eq:fdef}
    F(y)=-Y_+^{\tau+t(y)}y
  \end{equation}
and $t(y)\in(-\Delta,\Delta)$ is the unique time such that
\begin{equation}
  \label{eq:ftdef}
  \begin{split}
    \epsilon&=h\left(Y_-^{\,t(y)}y\right)\mbox{\quad if $h(y)\geq\epsilon$,}\\
    \epsilon&=h\left(Y_+^{\,t(y)}y\right)\mbox{\quad if $h(y)<\epsilon$.}    
  \end{split}
\end{equation}
\end{Theorem}
The expression~\eqref{eq:ftdef} of the traveling time $t(y)$ implies
that $F$ is continuous in $U(y_*(0))$ and smooth in each of its two
subdomains
\begin{equation}
  \begin{split}
    D_-&=U(y_*(0))\cap\{y:h(y)\geq\epsilon\} \mbox{\quad and}\\
    D_+&=U(y_*(0))\cap\{y:h(y)<\epsilon\}\label{eq:domains}
  \end{split}
\end{equation}
but, in general, its derivative
has a discontinuity along the boundary $D_0$ between $D_-$ and $D_+$.
The regularity of the implicit expression \eqref{eq:ftdef} for $t(y)$
for all $y\in U(y_*(0))$ follows from the strict transversality
\eqref{eq:stran:qg0} of the colliding periodic orbit $L$. Let us
denote the two smooth parts of the map $F$ by $F_+$ and $F_-$:
  \begin{equation}\label{eq:ftaupm}
    \begin{split}
      F_+(y)&=-Y_+^{\tau+t_+(y)}y\mbox{\quad where\quad}
      h\left(Y_+^{t_+(y)}y\right)=\epsilon\mbox{, and}\\
      F_-(y)&=-Y_+^{\tau+t_-(y)}y\mbox{\quad  where\quad}
      h\left(Y_-^{t_-(y)}y\right)=\epsilon\mbox{.}
    \end{split}
\end{equation}
That is, $F\vert_{D_+}=F_+$ and $F\vert_{D_-}=F_-$. Due to the regularity
of the definition of $t_\pm(y)$ in $U(y_*(0))$ the maps $F_\pm$ can
both be extended to the whole domain $D$. The map $F_+$ projects $D$
nonlinearly onto the local submanifold (the \emph{delayed switching manifold})
\begin{displaymath}
  \rg F_+=-Y_+^\tau[\{h=\epsilon\}\cap U(y_*(0))]
  =\{y\in U(y_*(0)):h(Y_+^{-\tau}[-y])=\epsilon\}\mbox{,}
\end{displaymath}
which has co-dimension one. The linearizations of $F_\pm$ with respect to
$y$ and the delay $\tau$ in $y=y_*(0)$ are
\begin{equation}\label{eq:linmap}
  \begin{split}
    F_\pm(y_*(0)+\eta;\tau+\theta)&=y_*(0)-A(\tau)
    \left[\left[\id-\frac{f_2 h_0'}{g_\pm}\right]\eta+\theta
      f_2\right] +O(|(\eta,\theta)|^2)\\
    g_-&=h_0'f_1=-h_0'A^\tau f_2\\
    g_+&=h_0'f_2
  \end{split}
\end{equation}
where $A(\tau)=\partial_y Y_+^\tau$ in $y_*(0)$, and
$f_1=f(y_*(0),-1,\lambda_0)$, $f_2=f(y_*(0),+1,\lambda_0)$, and
$h_0'=h'(y_*(0))$ (as introduced before). The linearization
of $F_+$ projects $\eta$ by the linear projection
$\id-f_2h_0'/(h_0'f_2)$ before propagating it with $A(\tau)$,
mirroring the dimension deficit of the image of the nonlinear map
$F_+$.

The map $F$ is a continuous piecewise smooth map in $\R^n$ with a rank
deficit in one half of the phase space. $F$ is implicitly defined and
nonlinear, even if the original relay system is piecewise linear of
the form \eqref{eq:lin}. Thus, Theorem~\ref{thm:retmapform} reduces
the study of the dynamics near a colliding symmetric periodic orbit to
the study of a low-dimensional piecewise smooth map in a similar
fashion as for the generic case \cite{S06}.  General bifurcation
theory has been developed for piecewise affine maps in $\R^n$, which
carries over partially to the nonlinear case
\cite{BG99,BBCK07,DFHH99,KDCHHPKN06,NOY94,SM07,SG06}.

\section{Single-degree-of-freedom oscillators}
\label{sec:osci}
In this section we demonstrate the use of the reduced maps derived in
Theorem~\ref{thm:retmapform}. We study the linear
single-degree-of-freedom oscillator \eqref{eq:osci} subject to a
delayed linear switch with hysteresis \eqref{eq:osci:disc}. We rescale
time $t$ and the variable $x$ in \eqref{eq:osci} such that the
equilibria of the flows $Y_\pm$ are at $\pm1$ and such that each flow
rotates with frequency $1$. Furthermore, we introduce the parameter
$\alpha$ which is the tilting angle of the switching decision function $h$.
This reduces system~\eqref{eq:osci} to a system with four parameters,
the damping $\zeta$, the delay $\tau\in(0,\infty)$, the width
$2\epsilon\in(0,\infty)$ of the relay hysteresis region, and the angle
$\alpha\in[-\pi/2,\pi/2]$ (negative feedback) of the normal vector to
the switching lines:
\begin{equation}
  \label{eq:oscired}
  \begin{split}
    \ddot x&+2\zeta\dot x+(1+\zeta^2) x=(1+\zeta^2)u\\
    u(t)&=  
    \begin{cases}
    -1 &\mbox{if $x(t-\tau)\cos\alpha+
      \dot x(t-\tau)\sin\alpha\geq\epsilon$}\\
    \phantom{-}1  &\mbox{if $x(t-\tau)\cos\alpha+
      \dot x(t-\tau)\sin\alpha\leq-\epsilon$}\\
    \lim_{s \nearrow t} u(s) &\mbox{if $x(t-\tau)\cos\alpha
      +\dot x(t-\tau)\sin\alpha\in(-\epsilon,\epsilon)$.}
  \end{cases}
  \end{split}
\end{equation}
We consider the case of an unstable focus (spiraling source)
corresponding to $\zeta<0$. In our numerical investigations we fix the
damping $\zeta$ to $-0.1$, which is a moderately expanding unstable
spiral. We study the dynamics near the colliding periodic orbit and
how it depends on the parameters of the control switch $u$ (delay
$\tau$, angle $\alpha$ and hysteresis width $\epsilon$).
\begin{figure}[t]
  \centering
  \includegraphics[scale=0.45]{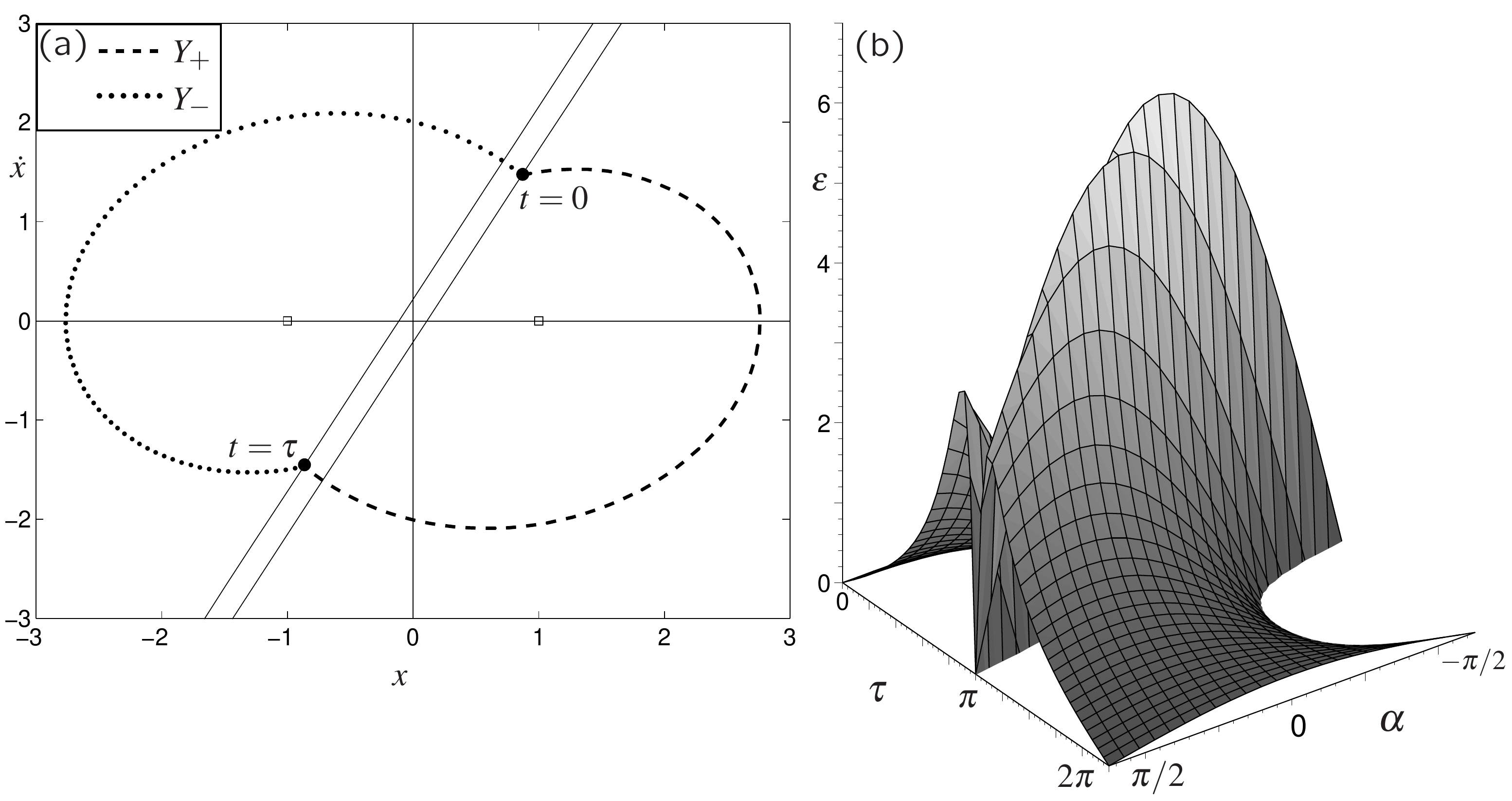}
  \caption{Panel (a): visualization of a colliding periodic orbit of
    \eqref{eq:osci},\,\eqref{eq:osci:disc}. Panel (b): collision
    surface in the three-dimensional parameter space
    $(\tau,\alpha,\epsilon)$.  Damping $\zeta$ is fixed at $-0.1$. At
    the parameters on the surface there exists a symmetric periodic
    orbit $L$ experiencing a corner collision. Beneath the surface the
    symmetric periodic orbit is a fixed point of the map $F_+$ in the
    subdomain $D_+$, above the surface it is a fixed point of the map
    $F_-$ in the subdomain $D_-$.}
  \label{fig:collsurf}
\end{figure}
The affine flows $Y_\pm$ are given by
\begin{displaymath}
  Y_\pm^ty=A(t)y\pm v(t)
\end{displaymath}
where
\begin{displaymath}
  \begin{split}
    A(t)&=e^{-\zeta t}\cos t-e^{-\zeta t}\sin t
    \begin{bmatrix}
      -\zeta&-1\\ 1+\zeta^2&\phantom{-}\zeta
    \end{bmatrix}\mbox{, and}\\
    v(t)&=e^{-\zeta t}
    \begin{bmatrix}
      -\zeta\sin t-\cos t +e^{\zeta t}\\ (1+\zeta^2)\sin t
    \end{bmatrix}
  \end{split}
\end{displaymath}
and have equilibria at $(\pm1,0)^T$.
The condition for the existence
of a colliding symmetric periodic orbit is
\begin{equation}
  \label{eq:collcond}
  \begin{split}
    y_1\cos\alpha+y_2\sin\alpha&=\phantom{-}\epsilon
    \mbox{\quad for $\tau>\pi$\ }\\
    y_1\cos\alpha+y_2\sin\alpha&=-\epsilon\mbox{\quad for $\tau<\pi$}
  \end{split}
\end{equation}
where $y=-[\id+A(\tau)]^{-1}v(\tau)$ (thus, $y=-Y_+^\tau y$).  For
$\zeta<0$, condition~\eqref{eq:collcond} ensures that the switching
point $y$ lies on the switching manifold and that $f_2=f(y,1)$ points
out of the hysteresis region. See Figure~\ref{fig:collsurf}(a) for a
visualization of a symmetric periodic orbit satisfying the collision
condition \eqref{eq:collcond}.
Figure~\ref{fig:collsurf}(b) shows the surface of parameters
$(\tau,\alpha,\epsilon)$ where a colliding symmetric orbit $L$ exists
in system~\eqref{eq:oscired}.  Whenever one varies the system
parameters along a path intersecting the surface transversally $L$
undergoes a corner collision. Beneath the surface $L$ is a fixed point
of $F_+$ in $D_+$, above the surface $L$ is a fixed point of $F_-$ in
$D_-$.
\begin{figure}[t]
  \centering
  \includegraphics[scale=0.7]{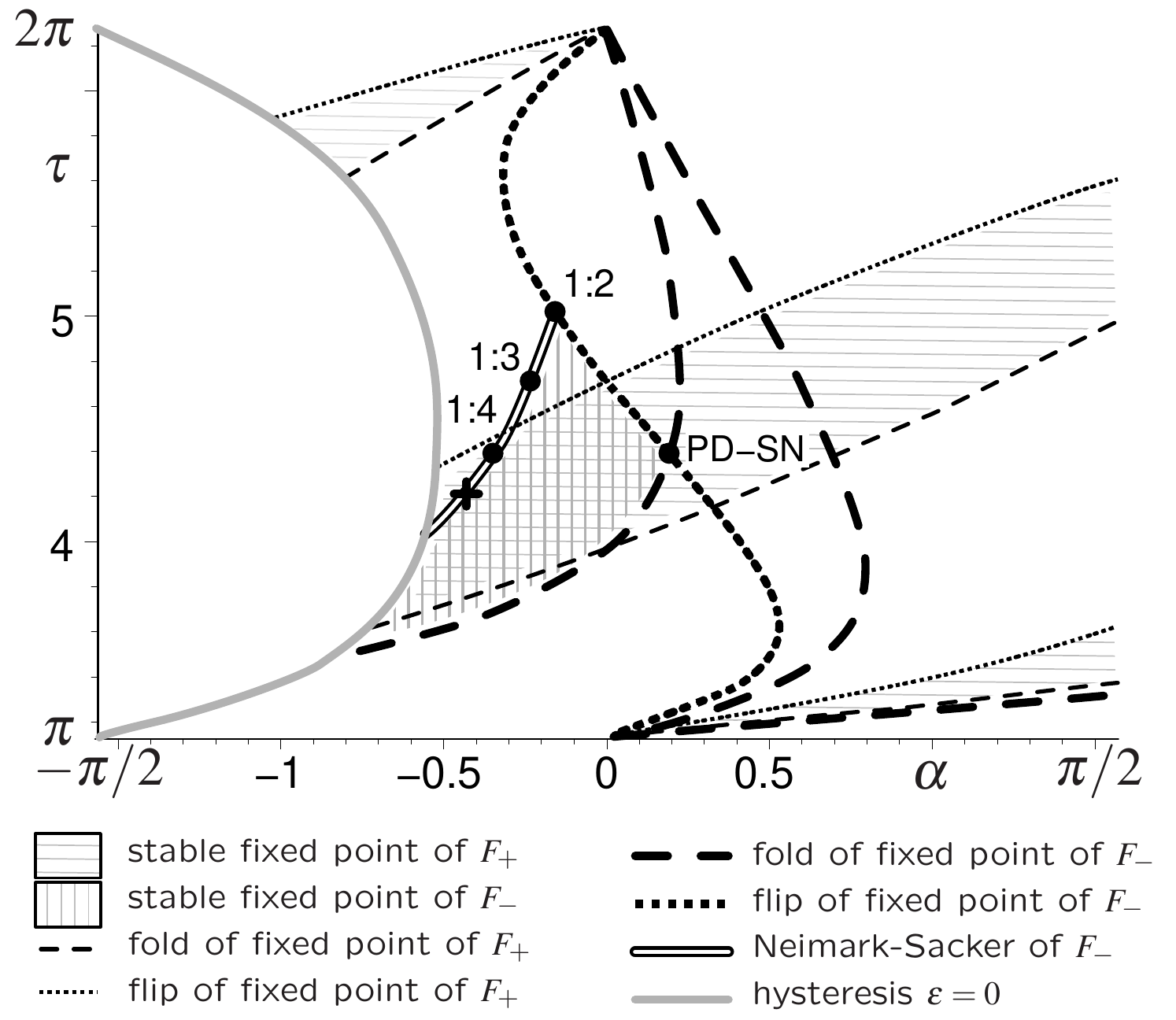}
  \caption{Bifurcations of the maps $F_+$ and $F_-$ projected onto the
    corner collision surface (`look from above' onto the surface in
    Figure~\ref{fig:collsurf}(b)).  Damping $\zeta$ is fixed at $-0.1$.
    Only the part of the collision surface with positive $\epsilon$ is
    shown. Thick lines are standard bifurcations of $F_-$. Thin lines
    are standard bifurcations of $F_+$. The cross is a coincidence of
    a Neimark-Sacker bifurcation of $F_-$ with the corner collision at
    $\epsilon=0.1$. Its numerical unfolding is presented in
    Figure~\ref{fig:bif2d} and Figure~\ref{fig:bif1d}.}
  \label{fig:locbif}
\end{figure}
The change of stability of the fixed point of $F$ at the corner
collision is determined by the linearizations of $F_-$ and $F_+$.
Figure~\ref{fig:locbif} shows a map of the different eigenvalue
configurations that can occur on the collision parameter surface for
$\tau\in(\pi,2\pi)$. The map $F_+$ has a one-dimensional image. Thus,
the linearization of $F_+$ in its fixed point has only one non-zero
eigenvalue. The horizontally hatched region in Figure~\ref{fig:locbif}
shows where this eigenvalue has modulus less than one. Within this
region the fixed point of $F_+$ is stable at collision. Consequently,
the symmetric periodic orbit of \eqref{eq:oscired} is linearly stable
for parameters beneath the horizontally hatched region in
Figure~\ref{fig:locbif} of the collision surface in
Figure~\ref{fig:collsurf}(b). The stable region of $F_+$ is bounded by a
flip bifurcation (eigenvalue equals $-1$, thin dotted in
Figure~\ref{fig:locbif}) and a fold bifurcation (or saddle-node,
eigenvalue equals $1$, thin dashed in Figure~\ref{fig:locbif}).
Parameter values on a bifurcation curve in Figure~\ref{fig:locbif}
correspond to codimension two events of the oscillator
\eqref{eq:oscired} because they also lie on the collision surface.
This means that the symmetric periodic orbit of \eqref{eq:oscired}, as a fixed
point of $F$, has a linearization with neutral stability in one of the
subdomains $D_\pm$ and, simultaneously, it is located on the boundary
between $D_-$ and $D_+$.

The map $F_-$ has a two-dimensional image. Thus, the linearization in
its fixed point has two potentially non-zero eigenvalues. The
vertically hatched region in Figure~\ref{fig:locbif} shows parameter
values where both eigenvalues are inside the unit circle. In this
region above the collision surface the oscillator has a symmetric
periodic orbit that is linearly stable as a fixed point of $F_-$. The
standard bifurcations of the fixed point of $F_-$ are shown as thick
lines (flip dotted, fold dashed, Neimark-Sacker hollow). Along the
Neimark-Sacker bifurcation (also called torus bifurcation) curve a
complex conjugate pair of eigenvalues is on the unit circle. All of
the bifurcation curves of $F_-$ correspond to codimension two events
for the symmetric periodic orbit of the oscillator \eqref{eq:oscired}
because they occur simultaneously with the collision, lying on the
collision surface in the three-dimensional parameter space shown in
Figure~\ref{fig:collsurf}(b).

\paragraph{Remarks} 
The flip bifurcation of $F_\pm$ corresponds to a symmetry breaking
bifurcation of the original symmetric periodic orbit of the oscillator
\eqref{eq:oscired} because the return map along the full periodic
orbit is the second iterate of $F$.

The line $\alpha=0$ in the figures \ref{fig:collsurf}(b) and
\ref{fig:locbif} corresponds to the case of pure position feedback
studied in \cite{CDHK06}. The linearizations of $F_+$ and $F_-$
coincide for $\alpha=0$. More precisely, the expressions for $g_+$ and
$g_-$ in \eqref{eq:linmap} are identical. It has been observed in
\cite{CDHK06} that non-smooth phenomena cannot occur for symmetric
periodic orbits at the event collision.

The points $(\alpha,\tau)=(0,\pi)$ and $(\alpha,\tau)=(0,2\pi)$ in
Figure~\ref{fig:locbif} are highly degenerate. The periodic orbit does
not intersect the switching line transversally at these parameter
values, violating the strict transversality
condition~\eqref{eq:stran:qg0}. The collision surface is singular for
$\tau=\pi$.

Figure~\ref{fig:locbif} shows codimension two degeneracies of the
linearization of $F_-$ such as a concurrence of eigenvalues $-1$ and
$+1$ at PD-SN and strong resonances along the Neimark-Sacker
bifurcation (double eigenvalue $-1$ at 1:2, eigenvalues $\exp(\pm2\pi
i/3)$ at 1:3, eigenvalues $\pm i$ at 1:4, see \cite{K04} for an analysis
and description). These points correspond to codimension three
bifurcations of the symmetric periodic orbit. Similarly, all crossings
of bifurcations of $F_-$ and bifurcations of $F_+$ in
Figure~\ref{fig:locbif} correspond to codimension three bifurcations
of the symmetric periodic orbit. Other bifurcations of higher
codimension that can occur along the bifurcation curves shown in
Figure~\ref{fig:locbif} involve the degeneracy of higher order terms
in the normal form. These special points have been omitted from
Figure~\ref{fig:locbif}.

The interaction PD-SN between the flip and the fold of the fixed point
of the map $F_-$appears degenerate in Figure~\ref{fig:locbif} in the
following sense. The flip curve crosses the fold curve transversally
instead of touching it quadratically as one would expect in a generic
parameter unfolding \cite{KMV04}. This is due to the projection of the
bifurcation curves onto the collision surface. The set of all
parameters in the $(\tau,\alpha,\epsilon)$-space where the fixed point
of $F_-$ is a saddle-node forms a surface. This surface does not
intersect the collision surface shown in Figure~\ref{fig:collsurf}(b)
transversally but touches it tangentially in the (dashed thick) fold
bifurcation curves of Figure~\ref{fig:locbif}. The same applies to the
fold of the fixed point of $F_+$.

\section{Unfolding of Neimark-Sacker bifurcation collision}
\label{sec:tor}
The symmetric periodic orbit $L$ of the oscillator \eqref{eq:oscired}
is stable near its corner collision in the hatched regions in
Figure~\ref{fig:locbif}. Where the two hatched regions overlap $L$ is
stable for parameters on both sides of the collision surface in
Figure~\ref{fig:collsurf}(b). Of primary interest are the parameter
regions near bifurcation curves bounding the region of stability of
$L$. In these regions we can expect that other (possibly stable) invariant
objects exist near $L$, which in turn collide with the boundary $D_0$
between the two subdomains $D_-$ and $D_+$ of the phase space. In
order to understand the dynamics near one of the codimension two
events one has to unfold it using two parameters.

A systematic classification of possible unfoldings of codimension two
bifurcations of piecewise smooth maps is not available in contrast to
the situation for smooth systems \cite{K04}. Due to the impossibility
of a general center manifold reduction it is difficult to derive
general results. The references \cite{BG99,BRG00,BBCK07,KDCHHPKN06} give a
long list of possible cases but discuss unfoldings only for very few
concrete examples.

\begin{figure}[t]
  \centering
  \includegraphics[scale=0.6]{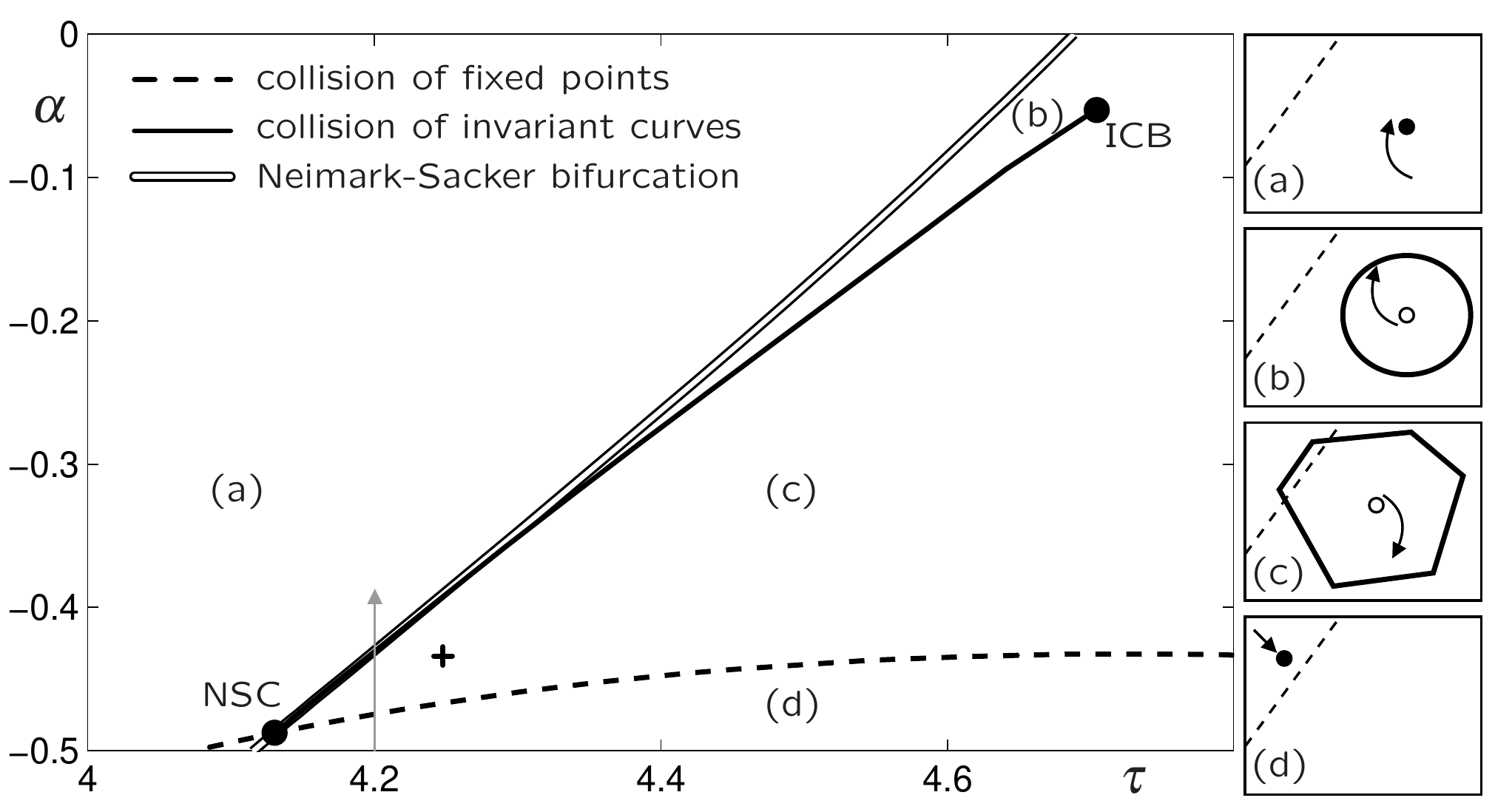}
  \caption{Unfolding of colliding Neimark-Sacker bifurcation in
    $(\tau,\alpha)$-plane.  Values of other parameters:
    $\epsilon=0.1$, $\zeta=-0.1$. All curves have been obtained by
    direct numerical continuation. The sketches on the right show
    qualitatively the dynamics in the different regions (a)--(d) (the
    dashed line is the switching line, attractors are solid, repellors
    hollow).  Figure~\ref{fig:starplots} shows the attracting
    invariant polygon at the parameter value marked by a cross.
    Figure~\ref{fig:bif1d} shows a parameter sweep of simulations
    along the gray arrow.}
  \label{fig:bif2d}
\end{figure}
In this section we describe in detail the dynamics near the symmetric
periodic orbit $L$ near a collision of a Neimark-Sacker bifurcation of $F_-$
with the boundary $D_0$. This case is only possible due to the
increase of the dimension of the phase space from one (in $D_+$) to
two (in $D_-$). In this sense it is the most characteristic feature of
the symmetric corner collision of the oscillator \eqref{eq:oscired}. It
is also the most complex case due to the involvement of invariant curves.

We fix $\epsilon=0.1$. The plane
$\{(\tau,\alpha,\epsilon):\epsilon=0.1\}$ intersects the collision
surface of Figure~\ref{fig:collsurf}(b) along a curve (not shown in
Figure~\ref{fig:locbif}) and it intersects the Neimark-Sacker
bifurcation on the collision surface (thick hollow curve) in a point
NSC. This point NSC is marked by a cross in Figure~\ref{fig:locbif}
and corresponds to a set of parameters where the periodic orbit $L$ is
on the boundary $D_0$, the linearization of $F_-$ has a complex
conjugate pair of eigenvalues on the unit circle, and the
linearization of $F_+$ is stable. We unfold this codimension two event
using the parameter $\alpha$ (tilting angle of the switching line)
and the delay $\tau$.  Figure~\ref{fig:bif2d} shows the bifurcation
diagram in the $(\tau,\alpha)$-plane. The point NSC corresponds to the
cross mark in Figure~\ref{fig:locbif}. The dynamics in the different
regions is sketched in the insets next to the diagram in
Figure~\ref{fig:bif2d}.

Let us explain the dynamics in the different regions near NSC (as
sketched in the insets to the right of the diagram in
Figure~\ref{fig:bif2d}). In parameter region (a) $L$ is in region $D_-$. $L$ is a
fixed point of $F_-$ and it is linearly stable, its linearization
having a pair of complex conjugate eigenvalues inside the unit circle.
On the Neimark-Sacker bifurcation curve (hollow) this pair of
eigenvalues is exactly on the unit circle. Thus, $L$ changes its
stability here.  Moreover, the Neimark-Sacker bifurcation of $L$ does
not have a strong resonance (see Figure~\ref{fig:locbif} where the
strong resonances are marked).  

\begin{figure}[t]
  \centering
  \includegraphics[scale=0.6]{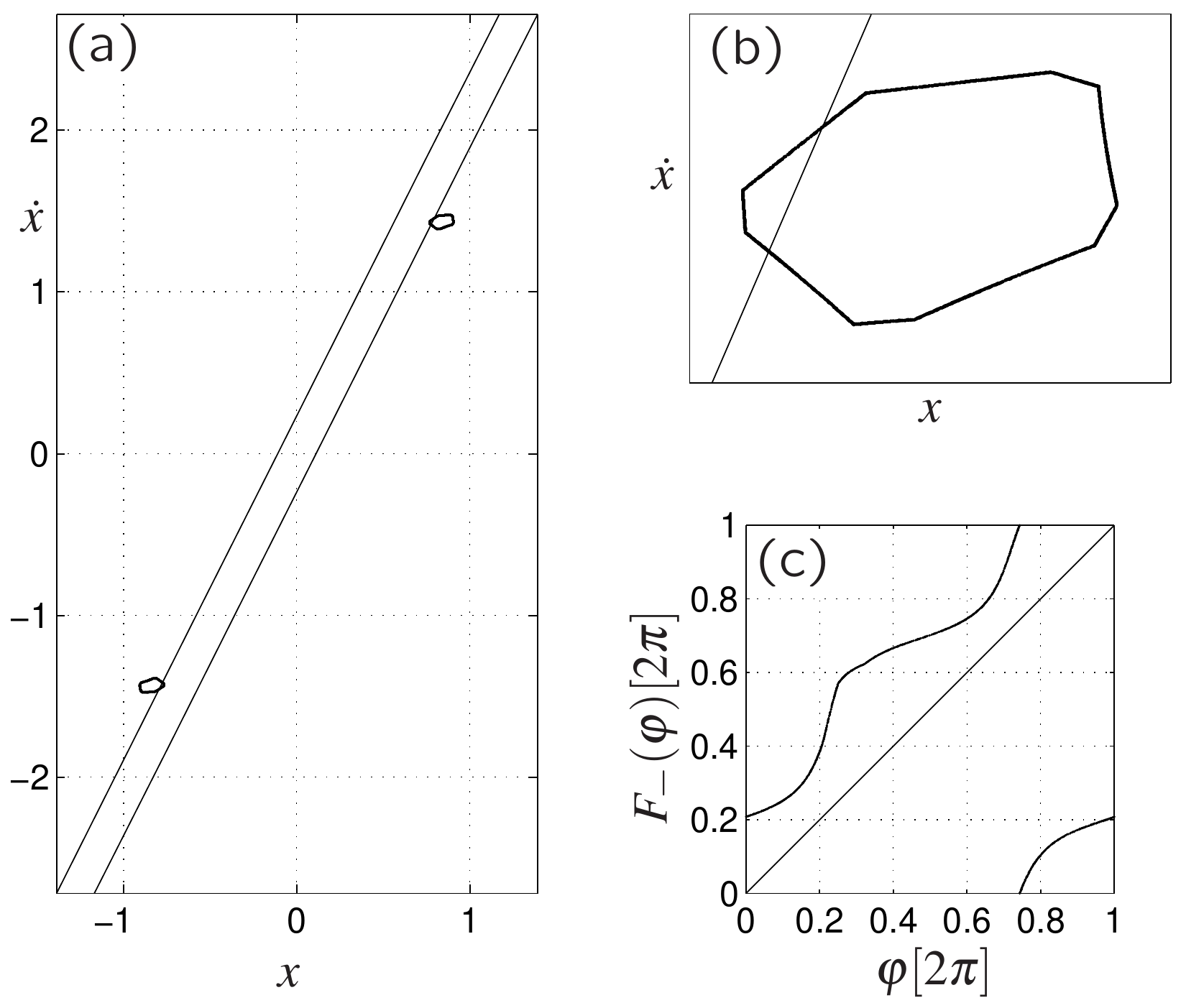}
  \caption{Attracting piecewise smooth invariant curve. Parameter
    values: $\alpha=-0.44$, $\tau=4.25$, $\epsilon=0.1$, $\zeta=-0.1$
    (cross in Figure~\ref{fig:bif2d}).  Panel (a) shows the switching
    points of the trajectory of the oscillator \eqref{eq:oscired} in
    the $(x,\dot x)$-plane. The switching lines are drawn as solid
    lines. Panel (b) zooms into the neighborhood of the attractor of $F$.
    Panel (c) shows the circle map on the closed polygon.}
  \label{fig:starplots}
\end{figure}
Consequently, a smooth closed invariant curve (invariant under $F_-$)
emerges from $L$ at the Neimark-Sacker bifurcation near NSC. In this
case the Neimark-Sacker bifurcation is supercritical, which implies
that the emerging closed invariant curve is stable and exists in the
region of linear instability of $L$ (region (b) in
Figure~\ref{fig:bif2d}). As the diameter of the invariant curve grows
it collides with the boundary $D_0$. This collision occurs along the
solid curve in Figure~\ref{fig:bif2d}.

At this collision the smooth closed invariant curve disappears. In
region (c) trajectories from all initial conditions close to $L$
(except $L$ itself) will eventually visit the region $D_+$. Thus, they
will eventually follow the map $F_+$ at least once getting projected
onto $\rg F_+$, a one-dimensional manifold. Since the fixed point of
$F_+$ is stable and in region $D_-$ (thus, it is not a fixed point of
$F$) all trajectories starting from points in $D_+$ will eventually
visit $D_-$. Hence, in region (c) all trajectories map back and forth
between the two regions $D_-$ and $D_+$.

This gives rise to invariant sets composed of finitely many smooth
arcs that are images of a section of $\rg F_+$ under $F_-$. Since
$F_-$ is approximately a rotation these arcs are rotations of $\rg
F_+$. Their composition forms a piecewise smooth closed invariant
curve consisting of finitely many arcs (a closed invariant polygon),
which is located partially in $D_-$ and partially in $D_+$. This type
of invariant set has been found and discussed for a piecewise linear
map already in \cite{SG06,SPG03}, and for flows in \cite{SO07,ZM06}.
Figure~\ref{fig:starplots} shows an example of such a polygon at the
parameter values marked by a cross in Figure~\ref{fig:bif2d}.
Figure~\ref{fig:starplots}(a) shows the switching points of the
trajectory of the original oscillator \eqref{eq:oscired} (where the
discrete component $u$ changes its value).
Figure~\ref{fig:starplots}(b) is a zoom into the vicinity of the upper
switching region, showing the attractor of $F$.
Figure~\ref{fig:starplots}(c) shows the map restricted to the
invariant polygon parametrized by the angle $\phi$ of the point on the
polygon relative to the the average of the attractor inside the curve.
The map is clearly invertible giving rise to either quasi-periodicity
or a pair of periodic orbits (locking) on the polygon. According to
\cite{SPG03} non-invertible maps on polygons (and, hence, in principle
chaotic dynamics) are also possible but we did not encounter this
phenomenon close to NSC. See also \cite{SG06,SPG03} for a study of
resonance and locking phenomena on the invariant polygons. 

If the parameters approach the collision curve of $L$ (dashed curve in
Figure~\ref{fig:bif2d}) the invariant polygons shrink as the fixed
points of $F_+$ and $F_-$ approach each other (and the boundary $D_0$).
Thus, in region (d) only the stable fixed point of $F_-$ in $D_-$
exists.

\begin{figure}[t]
  \centering
  \includegraphics[scale=0.6]{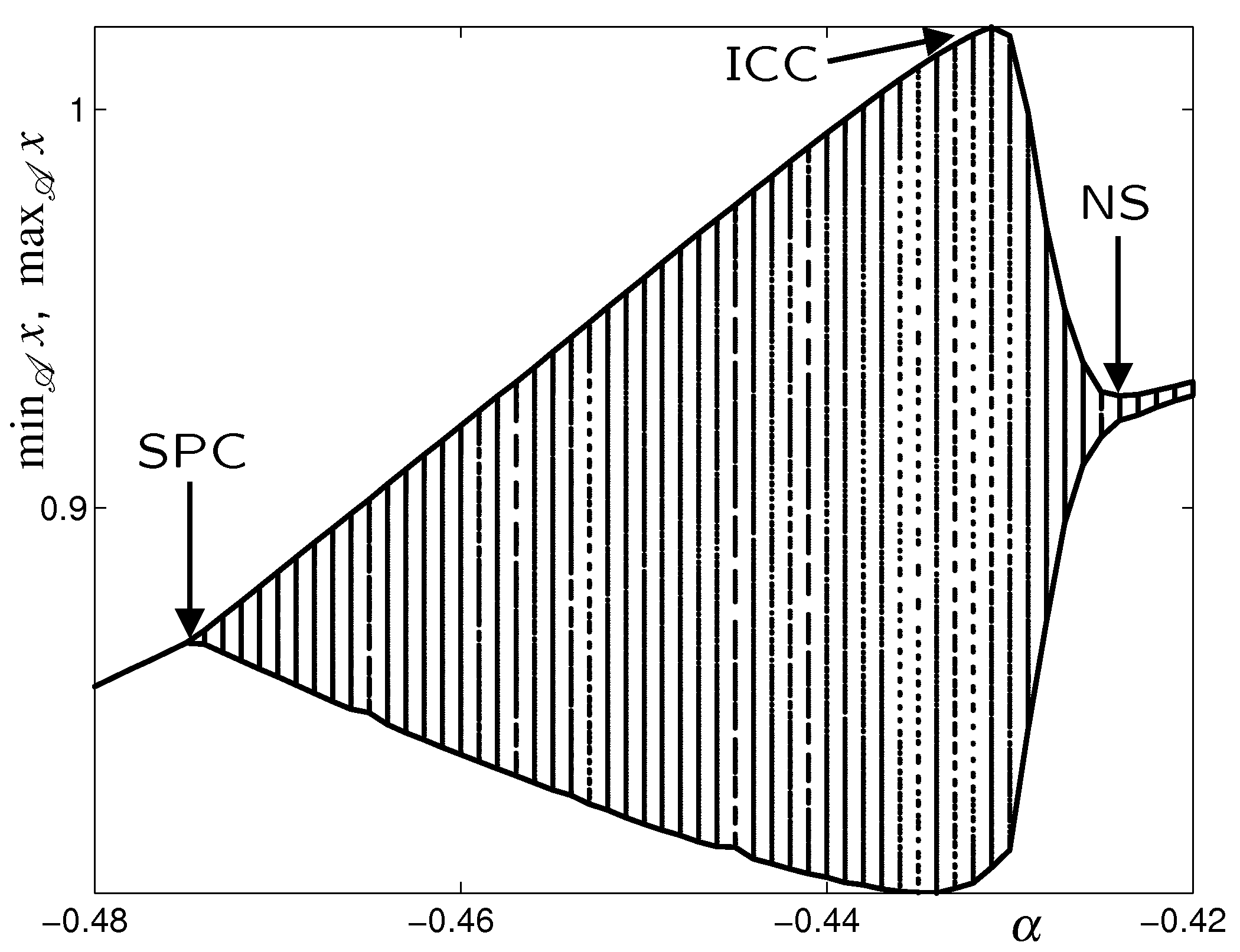}
  \caption{One-parameter pseudo bifurcation diagram for varying
    $\alpha$. Values of other parameters: $\tau=4.2$, $\epsilon=0.1$,
    $\zeta=-0.1$ (along the gray arrow of Figure~\ref{fig:bif2d}).
    Envelope contains maximum and minimum of the iterates of $F$
    between $40$ and $400$ iterations.  SPC, ICC and NS point out the
    intersection with the bifurcation curves in
    Figure~\ref{fig:bif2d}.}
  \label{fig:bif1d}
\end{figure}
Figure~\ref{fig:bif1d} shows the result of a parameter sweep along the
gray arrow in Figure~\ref{fig:bif2d} showing the maximum and the
minimum of the attractor ($\max F^j(y_0)$ and $\min F^j(y_0)$ for
$j=40,\ldots,400$, restarting such that $y_0$ is the last iterate from
the previous parameter). The figure gives evidence that the diameter
of the invariant polygons grows linearly in the parameter $\alpha$,
starting from the point SPC where the corner collision of $L$ occurs.
During the sweep the polygon also changes its shape (for example, the
number of corners or `overshoot' of its edges). A detailed analysis
how invariant polygons change under variation of the rotation and the
linear expansion of the two-dimensional map $F_-$ is carried out in
\cite{SO07}. The polygon attracts in finite time between the points
SPC and ICC in Figure~\ref{fig:bif1d} such that the numerical result
shown in the plot is very accurate for this region.  The detection of
the collision of the smooth invariant curve near ICC, the family of
smooth invariant curves and the Neimark-Sacker bifurcation NS are,
however, only inaccurately accessible by pure simulations due to the
weak attraction in region $D_-$.

\begin{figure}[t]
  \centering
  \includegraphics[scale=0.7]{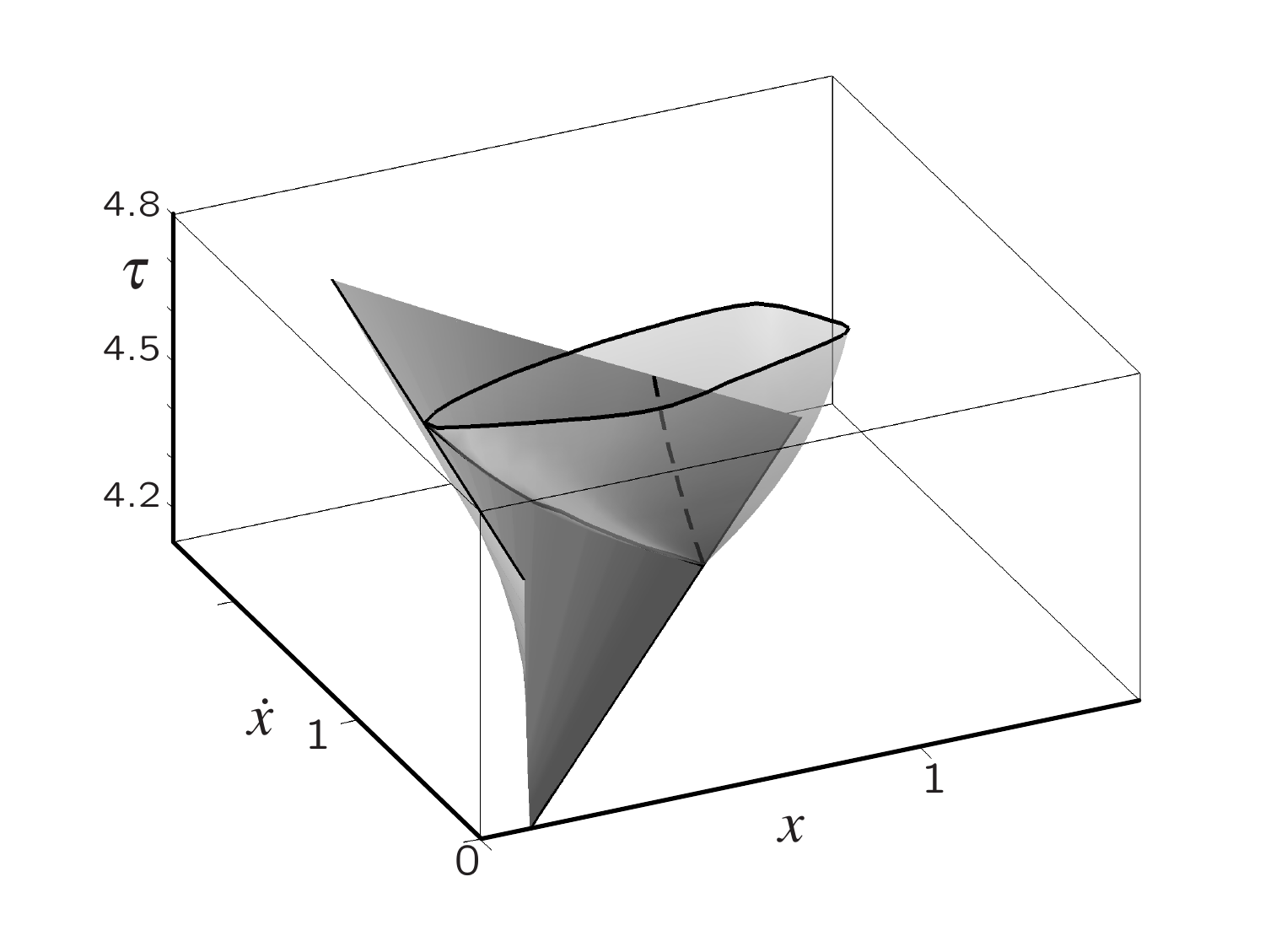}
  \caption{Family of closed invariant curves along collision curve
    (solid) in $(\tau,\alpha)$-plane of Figure~\ref{fig:bif2d},
    projected into the $(x,\dot x,\tau)$-space. The value of $\alpha$
    is coded by the gray value. Values of other parameters:
    $\epsilon=0.1$, $\zeta=-0.1$. The family of invariant curves is
    the transparent cone. The dashed curve inside the cone are the
    unstable equilibria of $F$. The fanned out intransparent surface
    are the switching lines. The points where the invariant curves
    touch the switching lines are marked as a solid black curve lying
    on the cone and the intransparent surface.}
  \label{fig:tc-family}
\end{figure}
Neimark-Sacker bifurcations and colliding smooth invariant curves can
be efficiently computed directly because only the smooth
two-dimensional map $F_-$ and the expression for the boundary $D_0$
are involved. The Neimark-Sacker bifurcation in Figure~\ref{fig:bif2d}
can be accurately and efficiently continued using standard numerical
algorithms as described in detail in \cite{K04} and implemented in
AUTO \cite{DCFKSW98}. The smooth closed invariant curve is given implicitly
by the invariance equation
\begin{equation}
  \label{eq:invcurve}
    y(\eta(\phi))=-Y_+^{\,\tau+t(\phi)}y(\phi)\mbox{,\quad}
    \epsilon=
      [\cos\alpha, \sin\alpha]
\cdot Y_-^{\,t(\phi)}y(\phi)
\end{equation}
where $y(\phi)$ is given by
\begin{equation}\label{eq:param}
    y(\phi)=y_0+r(\phi)
    \begin{bmatrix}
      \cos\phi\\\sin\phi
    \end{bmatrix}
\end{equation}
and $y_0$ is the equilibrium of $F_-$:
\begin{equation}
  \label{eq:nseq}
    y_0=-Y_+^{\,\tau+t_0}y_0\mbox{,\quad}
    \epsilon=[\cos\alpha,\sin\alpha]\cdot Y_-^{\,t_0}y_0 \mbox{.}
\end{equation}
The equations \eqref{eq:invcurve}--\eqref{eq:nseq} have the periodic
functions $r(\phi)$, $\eta(\phi)$ and $t(\phi)$, the vector $y_0\in\R^2$ and
the scalar $t_0$ as variables. The definition of $y(\phi)$ via
\eqref{eq:param} corresponds to a parametrization of the invariant
curve by angle with respect to the equilibrium $y_0$. The function $r$
is the distance of the point on the invariant curve from $y_0$, $\eta$
is the circle map and $t$ is the time elapsed from the last switch
(the map $F_-$ is defined only implicitly). This parametrization can
be expected to work only in two-dimensional maps and for convex
curves. Close to the point NSC in Figure~\ref{fig:bif2d} (in fact,
close to the Neimark-Sacker curve) the invariant curves are ellipses
and, thus, convex, making the parametrization \eqref{eq:param}
regular. We extend system~\eqref{eq:invcurve}--\eqref{eq:nseq} by 
the collision condition
\begin{equation}
  \label{eq:invcollcond}
  0=\max_{\phi\in(0,2\pi]}t(\phi)\mbox{.}
\end{equation}
The right-hand-side of \eqref{eq:invcollcond} is a smooth function of
$\phi$ close to the point NSC in Figure~\ref{fig:bif2d} due to the
convexity of the invariant curve. Thus, the system of equations
\eqref{eq:invcurve}--\eqref{eq:invcollcond} for the variables
$(r,\eta,t,y_0,t_0,\tau,\alpha)$ defines a smooth family of colliding
closed invariant curves which can be found by a Newton iteration with
pseudo-arclength embedding \cite{KASP85}. Its projection into the
$(\tau,\alpha)$-plane is shown as a solid curve in
Figure~\ref{fig:bif2d}. The closed invariant curves of the colliding
family projected into the space $(x,\dot x,\tau)$ are depicted (as a
transparent gray cone) in Figure~\ref{fig:tc-family}.
Figure~\ref{fig:tc-family} gives evidence of the linear dependence of
the radius $r$ of the closed invariant curve on the parameter $\tau$.
Thus, the solid collision curve must be quadratically tangent to the
Neimark-Sacker curve at NSC in Figure~\ref{fig:bif2d} because the
radius $r$ has square root like asymptotics with respect to the
distance of the parameter value to the Neimark-Sacker bifurcation. At
the point ICB the invariant curve is close to a break-up such that the
numerical approximation (32 complex Fourier modes) of $r$ and $\eta$
gives an error estimate greater than $10^{-2}$.

\paragraph*{Remarks}
If the Neimark-Sacker bifurcation is subcritical then the line of
collisions of closed invariant curves (solid) lies on the other side
of the Neimark-Sacker bifurcation (hollow) in Figure~\ref{fig:bif2d}
(but still tangential). The polygons remain stable filling
the region between the solid and the dashed curve.

The interaction between flip bifurcation and collision (dotted curves
in Figure~\ref{fig:locbif}) has been analyzed in \cite{KDCHHPKN06}. A
periodic orbit $L_2$ of period two branches off from $L$.  In any
generic two-parameter unfolding a curve of collisions of $L_2$ emerges
tangentially to the flip curve of $L$ (looking very similar to
Figure~\ref{fig:bif2d}, replacing the Neimark-Sacker curve by a flip
curve and the collision of the invariant curve by a collision of
$L_2$). The dynamics of the system between the collision of $L_2$ and
the collision of the period-one orbit $L$ has, for example in
\cite{KDCHHPKN06}, a stable period-two orbit flipping between the
subdomains $D_-$ and $D_+$.

\section{Conclusions and open problems}
\label{sec:conc}
The local bifurcation theory of periodic orbits in hybrid dynamical
systems with delayed switching can be reduced to the study of
low-dimensional maps because the local return map (Poincar{\'e} map)
of periodic orbits is finite-dimensional. 
The codimension one events for generic systems (corner
collision and tangential grazing) were studied in \cite{S06}. This
paper studies the case of simultaneous corner collision due to
reflection symmetry, a case that has been encountered frequently in
example studies \cite{BKW06,BH98,CDHK06,H91}.

The presence of hysteresis in the switch simplifies the proof of the
reduction theorem (Theorem~\ref{thm:retmapform}), given here for the
first time. However, the statement of the reduction theorem remains
valid also without hysteresis.

The derivation of Theorem~\ref{thm:retmapform} can be generalized in a
straightforward manner to other discrete symmetries and other than
binary switches (with more notational overhead) as long as the
symmetry can be reduced locally near the periodic orbit.

The reduction extends the applicability of theory and numerical
methods that have been developed for smooth and piecewise smooth
low-dimensional maps to systems with delayed switches.  On the
numerical side this includes direct continuation of periodic orbits
and their bifurcations and discontinuity induced events (such as
grazing and collision) and the continuation of smooth invariant
curves.  Robust and universal methods for continuation and detection
of discontinuity induced bifurcations for periodic orbits have been
developed by Piiroinen \cite{PVC04,PK07}. Methods for the continuation
of closed invariant curves have a longer history. Schilder
\cite{SOV05} and Thakur \cite{DT06} present recent implementations,
the papers also give surveys of earlier work. There is, however, a
large gap between recent developments of numerical methods for closed
invariant curves and piecewise smooth systems and their actual
availability in the form of software. Due to this gap the investigation
of the oscillator in Section~\ref{sec:tor} could not rely on
generally available tools. 

The analysis of the oscillator shows that dynamical phenomena of
hybrid systems with delayed switches can be systematically discovered
with the help of numerical continuation and the reduction theorem.
Two major open problems can be identified from the results of this
analysis and the comparison to other results.

First, there is a gap between the relatively simple local bifurcation
theory of periodic orbits as presented here and in \cite{S06} and the
abundance and variety of complex phenomena observed in detailed
example studies such as \cite{BKW06,CDHK06,H91}. In particular, many
of the complex phenomena (say, chaotic attractors that have periodic
orbits with various switching patterns embedded) cannot be reached by
systematic continuation of periodic orbits and branching off at
continuous local bifurcations. One of the reasons behind this gap is
the presence of discontinuous events such as the grazing and the
corner collision with $q<0$ presented in Section~\ref{sec:genbif}.
These events are beyond the scope of local bifurcation theory because
some trajectories leave the local neighborhood of the periodic orbit.
A way to close this gap may be the study of discontinuous events under
the assumption of the existence of a smooth `global' return map for
trajectories leaving the neighborhood similar to the treatment of
global bifurcations of smooth dynamical systems.

The second open question is the robustness of discontinuity induced
bifurcations with respect to singular perturbations. For example, do
the closed invariant polygons (as shown in Figure~\ref{fig:starplots})
persist when we change the oscillator~\eqref{eq:osci} to
\begin{equation}
  \label{eq:osciext}
  \ddot x+\zeta y +\omega^2 x=u\mbox{,\quad} \beta \dot y=\dot x-y
\end{equation}
(a modification to non-viscous damping for visco-elastic materials
\cite{A05}) with a small $\beta$? In particular, are there still
finitely many smooth arcs? Note that the argument of \cite{SPG03} and
of Section~\ref{sec:tor} that the arcs are images under $F_-$ of the
one-dimensional image $\rg F_+$ cannot be applied anymore because
$\dim\,\rg F_+=2$ for \eqref{eq:osciext}. Most statements of the
bifurcation theory for piecewise smooth systems cannot be easily
generalized to higher dimensions due to the lack of center manifolds.
For example, the limit of small delay $\tau$, which is governed by
standard singular perturbation theory for smooth dynamical systems
\cite{C03}, is highly non-regular for \eqref{eq:osci} in the case of
$\epsilon=0$, see \cite{S06}.

\section*{Acknowledgments}
We wish to thank our colleagues Petri Piiroinen and Robert Szalai for
fruitful discussions initiating this paper. The research of J.S. and
P.K. was partially supported by by EPSRC grant GR/R72020/01.
\bibliographystyle{plain}
\bibliography{delay} 
\begin{appendix}
\section{Basic properties of the forward evolution}
\label{sec:app:evolution}
Let $\xi_0\in C([-\Theta,0];\R^n)$ be the initial history segment of the
continuous variable $y$ (where $\Theta\geq\tau$) and $u_0\in\{-1,1\}$
be the initial state of the discrete variable $u$ in
system~\eqref{eq:cont},\,\eqref{eq:disc}. How does this initial state
evolve to time $T$ (defining the evolution $E^T(\xi_0,u_0)$?

First, we define $E^T(\xi_0,u_0)$ for $T\in(0,\tau]$
using the variation of constants formulation of \eqref{eq:cont}. We
define the following subsets of the closed interval $[-\tau,0]$:
\begin{displaymath}
  \begin{split}
    R_+&=\{s\in[-\tau,0]:h(\xi_0(s))\leq-\epsilon\}\\
    R_-&=\{s\in[-\tau,0]:h(\xi_0(s))\geq\epsilon\}\\
    R_{0\ }&=\{s\in[-\tau,0]:h(\xi_0(s))\in(-\epsilon,\epsilon)\}\mbox{.}
  \end{split}
\end{displaymath}
The set $R_0$ is open relative to $[-\tau,0]$. Thus, if non-empty it
is a union of countably many disjoint open (relative to $[-\tau,0]$) intervals.
We arrange this sequence of countably many disjoint intervals into two
subsequences of intervals: $I^+_j$ ($j=1,\ldots$), $I^-_j$
($j=1,\ldots$) and, possibly, one extra interval $I_0$:
\begin{displaymath}
  \begin{split}
    I^+_j&=\mbox{intervals of $R_0$ that have a lower boundary $s^+_j\in R_+$, 
      i.e., $h(\xi_0(s^+_j))=-\epsilon$,}\\
    I^-_j&=\mbox{intervals of $R_0$ that have a lower boundary $s^-_j\in R_-$, 
      i.e., $h(\xi_0(s^-_j))=\epsilon$,}\\
    I_0&=
    \begin{cases}
      [-\tau,t_0)\quad\mbox{
        \begin{minipage}[t]{0.7\linewidth}
          if $-\tau\in R_0$ (then $t_0$ is the upper boundary of the
          left-most interval of $R_0$)
        \end{minipage}
      }\\
      \emptyset \quad\mbox{if $-\tau\notin R_0$.}
    \end{cases}
  \end{split}
\end{displaymath}
We define the following function $\chi:[0,\tau]\mapsto\R$:
\begin{equation}\label{eq:chi}
  \chi(s)=
  \begin{cases}
    \phantom{-}1 &\quad\mbox{if $s-\tau\in R_+\cup \bigcup_jI^+_j$}\\
    -1 &\quad\mbox{if $s-\tau\in R_-\cup \bigcup_jI^-_j$}\\
    \phantom{-}u_0 &\quad\mbox{if $s-\tau\in I_0$.}
  \end{cases}
\end{equation}
Thus, $\chi$ is measurable on $[0,\tau]$ and either $1$ or $-1$ everywhere.
The variation-of-constants formulation of \eqref{eq:cont} is
\begin{equation}\label{eq:voc}
  y(t)=\xi_0(0)+\int_0^tf(y(s),+1)\frac{[1+\chi(s)]}{2}+f(y(s),-1)\frac{[1-\chi(s)]}{2}\,ds.
\end{equation}
This is a fixed point problem on the space $C([0,\tau];\R^n)$ of
continuous functions on the interval $[0,\tau]$ that has a globally
unique solution $y(\cdot)\in C([0,\tau];\R^n)$ due to the Lipschitz
continuity of $f(\cdot,\pm1)$. Then, for $T\in(0,\tau]$,
$E^T(\xi_0,u_0)=(\xi(T)(\cdot),u(T))$ is defined by
\begin{equation}\label{eq:evol}
  \begin{split}
    \xi(T)(s)&=
    \begin{cases}
      y(T+s) &\quad\mbox{if $s\in[-T,0]$,}\\
      y_0(T+s) &\quad\mbox{if $s\in[-\tau,-T]$,}
    \end{cases}\\
    u(T)&=\chi(T)\mbox{.}
  \end{split}
\end{equation}
We observe that the initial value of the discrete variable $u_0$ only
affects the result $E^T(y_0,u_0)$ if
$h(\xi_0(-\tau))\in(-\epsilon,\epsilon)$. Otherwise, \eqref{eq:evol}
simply sets the discrete variable to its consistent value. For $T>\tau$
we define $E^T(\xi_0,u_0)$ as a concatenation of smaller time steps, for
example, if $T\in[(k-1)\tau,k\tau]$ then $E^T=E^{T/k}\circ
\ldots \circ E^{T/k}$. This definition is independent of the
particular partition of $[0,T]$.

The definition of the evolution $E$ using the variation-of-constants
formulation \eqref{eq:voc} allows one to initialize $E$ from arbitrary
continuous history segments $\xi_0$ and discrete states $u_0$ even if
$\xi_0$ crosses the switching manifolds $\{h(x)=\pm1\}$ infinitely
often or if $u_0$ is `inconsistent'. Due to \eqref{eq:voc} the
continuous component $E_c$ of $E^T(y_0,u_0)$ depends continuously on
$T$. If $T\geq\tau$ then $E_c^t(\xi_0,u_0)$ is Lipschitz
continuous with respect to $t$. Its Lipschitz constant is
$L_{\max}\cdot\|E_c^t(\xi_0,u_0)\vert_{t\in[t_0,t_E]}\|$ where $L_{\max}$ is the
Lipschitz constant of the right-hand-side $f$. For positive times the
trajectory of the discrete component $E_d$ is continuous from the
right, that is $\lim_{s\searrow t}E_d^s(\xi_0,u_0)=E_d^t(\xi_0,u_0)$
for all $t>0$.

\paragraph*{Proof of Lemma~\ref{thm:basic}}
(\emph{Point~\ref{thm:swfin}})\qquad Let $\delta>0$ be such that
$|s-t|<\delta$ implies $|h(\xi_0(s))-h(\xi_0(t))|<2\epsilon$ for all
$s,t\in[-\tau,0]$. This $\delta$ exists because the initial value
$\xi_0\in C([-\tau,0];\R^n)$ is uniformly continuous and $h$ is
Lipschitz continuous.

We have to check how many sign changes the function $\chi$, defined in
\eqref{eq:chi}, can have in the interval $[0,\tau]$. Let us denote the
upper boundary of each interval $I^\pm_j$ by $t^\pm_j$. If $I_0$ is
non-empty $\chi$ can change its sign only in $[t_0,\tau]$ (one change
is, possibly, in $t_0$). If $I_0$ is empty we will use the notation
$t_0=0$ in the following argument.

After $t_0$ the function $\chi$ changes from $-1$ to $1$ only at times
$t=t^-_j+\tau$ when $t^-_j\in R_+$, that is, $h(\xi_0(s^-_j))=\epsilon$
(by definition of $I^-_j$) and $h(\xi_0(t^-_j))=-\epsilon$. By
definition of $\delta$ this implies that $t^-_j\geq s^-_j+\delta$.
Hence, $\chi=-1$ for at least time $\delta$ before it can switch to
$1$. The same argument applies for intervals $I^+_j$ and $\chi$
switching from $1$ to $-1$.  Consequently, between two subsequent
switchings of $\chi$ a time of at least $\delta$ must elapse. This
limits the number of switchings to a finite number on a bounded
interval.

\emph{Point 2}: After time $\tau$ the constant $\delta$ in the above
argument is bounded from below by
\begin{equation}\label{eq:mindist}
  \delta\geq\frac{2\epsilon}{H_{\max}L_{\max}\max\|y\|}
\end{equation}
where $y$ is the solution of the fixed point problem for the
variation-of-constants formulation \eqref{eq:voc}, $H_{\max}$ is
the Lipschitz constant of switching function $h$, and $L_{\max}$ is
the Lipschitz of the right-hand-side $f$.

\paragraph*{Proof of Lemma~\ref{thm:continuity} (continuity)}
It is sufficient to prove the continuity of $E_c^T(\xi,u)$ with
respect to $\xi$ in $\xi_*=y_*(t+\cdot)\in C([-\Theta,0];\R^n)$ and
$u=u_*(t)$ for times $T\leq\tau$ because $(y_*,u_*)$ is periodic. 

Due to the Lipschitz continuity of $f(\cdot,\pm1)$ in the
variation-of-constants formulation \eqref{eq:voc} it is sufficient to
prove that, for any given $\Delta>0$, we can find a neighborhood
$U(\xi_*)$ such that
\begin{equation}\label{eq:discdiff}
  \int_0^T|E_d^s(\xi,u_*(t))-u_*(t+s)|\,ds<2\Delta
\end{equation}
for all $\xi\in U(\xi_*)$. That is, we have to show that, starting
from $\xi$, we follow the same flow as $y_*$ all the time in $[0,T]$
except in a union of intervals of overall length $\Delta$. Let
$t_{*,k}$ ($k=1,\ldots,\mu$) be the crossing times of $\xi_*$ in
$[-\tau,0]$. Let us denote by $u_k$ the value of the discrete variable
$u_*$ at time $t$ and at the switching times in $[t,t+\tau]$, that is,
$u_0=u_*(t)$, $u_k=u_*(t+\tau+t_{*,k})$ ($k=1,\ldots,\mu$). Since
$u_*$ is continuous from the right $u_*(t+\tau+s)=u_k$ also for $s$
slightly larger than $-\tau$ and $t_{*,k}$ ($k=1,\ldots,\mu$).

Due to the weak transversality of $(y_*,u_*)$ for all sufficiently
small $\delta_1>0$ there exists a $\delta_2>0$ such that
\begin{enumerate}
\item\label{prf:bigger} $-\delta_2\geq u_k\cdot h(y_*(t+t_{*,k}+\delta_1))+\epsilon$
  for all $k=1,\ldots,\mu$, and
\item\label{prf:smaller} $-\delta_2\geq u_{k-1} h(y_*(t+s))-\epsilon$
  for all $s\in[t_{*,k-1},t_{*,k}-\delta_1]$ (for $k=2,\ldots,\mu$)
  and $s\in[-\tau,t_{*,1}-\delta_1]$.
\end{enumerate}
Both statements follow from the strict monotonicity of
$|h(y_*(\cdot)|$ at crossing times of $(y_*,u_*)$ and the identities
$u_kh(y_*(t+t_{*,k}))+\epsilon=0$ and
$u_{k-1}h(y_*(t+t_{*,k}))-\epsilon=0$.  The points \ref{prf:bigger}
and \ref{prf:smaller} imply that for all sufficiently small
$\delta_1>0$ there exists an open neighborhood $U(\xi_*)$ such that
$\xi\in U(\xi_*)$ satisfy
\begin{enumerate}
\item\label{prf:bigger2}  $0\geq u_k\cdot h(\xi(t_{*,k}+\delta_1))+\epsilon$
  for all $k=1,\ldots,\mu$, and
\item\label{prf:smaller2} $0\geq u_{k-1} h(\xi(s))-\epsilon$ for all
  $s\in[t_{*,k-1},t_{*,k}-\delta_1]$ (for $k=2,\ldots,\mu$) and
  $s\in[-\tau,t_{*,1}-\delta_1]$.
\end{enumerate}
Consequently, for initial values $(\xi,u)$ with $\xi\in U(\xi_*)$ and
$u=u_*(t)$ the discrete variable $E_d^s(\xi,u)$ is equal to $u_0$ in
$[0,\tau+t_{*,1}-\delta_1]$ and equal to $u_k$ in
$[\tau+t_{*,k}+\delta_1,\tau+t_{*,k+1}-\delta_1]$ for
$k=1,\ldots,\mu-1$ and, if $t_{*,\mu}+\delta_1<0$, equal to $u_\mu$ in
$[\tau+t_{*,\mu}+\delta_1,\tau]$. Thus, $E_d^s(\xi,u)$ and $u_*(t+s)$
can be different only in $\mu+1$ intervals of length less than
$2\delta_1$. Choosing $\delta_1=\Delta/(\mu+1)$ we obtain the neighborhood
$U(\xi_*)$ guaranteeing \eqref{eq:discdiff} for all $\xi\in U(\xi_*)$. 

\paragraph*{Proof of Lemma~\ref{thm:finmap}}
First we choose a $\delta>0$ sufficiently small such that $|h(y_*(s))|$
is strictly monotone increasing in all intervals
$[t_{*,k}-\delta,t_{*,k}+\delta]$ where $t_{*,k}$ ($k=1,\ldots,\mu$) are
the crossing times of the periodic orbit $(y_*,u_*)$ in $(-\tau,0)$. This $\delta$
exists due to the weak transversality of $(y_*,u_*)$. We define the constant
\begin{displaymath}
  c_0=\min_{k=1,\ldots,\mu}\left\{|h(y(t_{*,k}-\delta))|-\epsilon|,
    |h(y(t_{*,k}+\delta))|-\epsilon|\right\}\mbox{,}
\end{displaymath}
which is positive. Thus, $|h(y_*(t_{*,k}-\delta)|-\epsilon\leq-c_0$
and $|h(y_*(t_{*,k}+\delta))|-\epsilon\geq c_0$.
We choose $U_1(\xi_*)$ such that for all $\xi\in U_1(\xi_*)$
\begin{displaymath}
  \max_{k=1,\ldots,\mu}\left\{|h(\xi(t_{*,k}-\delta))-h(y_*(t_{*,k}-\delta))|,
    |h(\xi(t_{*,k}+\delta))-h(y(t_{*,k}+\delta))|\right\}\leq c_0/2\mbox{.}  
\end{displaymath}
This guarantees that $|h(\xi(\cdot))|-\epsilon$ changes its sign
in $[t_{*,k}-\delta,t_{*,k}+\delta]$ for all $\xi\in U(\xi_*)$.
Consequently, there exist unique times $t_k$ ($k=1,\ldots,\mu$) such
that
\begin{displaymath}
  t_k=\min\{t\in[t_{*,k}-\delta,t_{*,k}+\delta]:|h(\xi(t))|=\epsilon\}
  \mbox{,\quad($k=1,\ldots,\mu$).}
\end{displaymath}
If $\delta$ is sufficiently small $E_d^t(\xi,1)$ will change its value
exactly at the times $t_k+\tau$ in $[0,\tau]$. This implies that
$\xi^\tau=E_c^\tau(\xi,1)$ is given by the recursion (denoting
$t_0=-\tau$ and $t_{\mu+1}=0$)
\begin{displaymath}
  \begin{split}
    \xi^\tau(-\tau)&=\xi(0)\\
    \xi^\tau(t)&=
    \begin{cases}
      Y_-^{[t-t_k]}\xi^\tau(t_k) &\mbox{if $t\in(t_k,t_{k+1}]$ and $k$ odd}\\
      Y_+^{[t-t_k]}\xi^\tau(t_k) &\mbox{if $t\in(t_k,t_{k+1}]$ and $k$ even}\mbox{.}
    \end{cases}
  \end{split}
\end{displaymath}
Hence, $\xi^\tau=E_c^\tau(\xi,1)$ depends only on $\xi(0)=y(0)$ and
$t_1,\ldots,t_\mu$. The discrete variable $E_d^\tau(\xi,1)$ equals
$u_*(\tau)$. Therefore, $E^p(\xi,1)=E^{p-\tau}\circ E^\tau(\xi,1)$ also
depends only on $\xi(0),t_1,\ldots,t_\mu$. 

Denote the continuous component $E_c^p(\xi,1)$ of $E^p(\xi,1)$ by
$\xi^p$. The cross section $\Sigma$ is transversal to $y_*(\cdot)$ at
$y_*(0)$.  Thus, for every $y\in U(y_*(0))$ there exists a locally
unique traveling time $t$ such that $Y_+^ty\in\Sigma$. The function
$t_\Sigma:y\mapsto t$ is well defined (and smooth) in $U(y_*(0))$. If
$U_1(\xi_*)$ is sufficiently small then $\xi^p(0)\in U(y_*(0))$ due to
the continuity of $E_c$ in $\xi_*=y_*(\cdot)$. The return time $T(\xi)$ is given
by $T(\xi)=p+t_\Sigma(\xi^p(0))$. As $\xi^p$ depends only on
$(\xi(0),t_1,\ldots,t_\mu)$ the same applies to $T(\xi)$ and, hence,
the return map $P$.

\paragraph*{Proof of Corollary~\ref{thm:smoothmap}}
The condition \ref{cond:nocoll} of Corollary~\ref{thm:smoothmap}
guarantees that the periodic orbit follows one of the flows in each of
its crossing times $t_{*,k}\in(0,p)$ ($k=1,\ldots,m$), say,
$Y_k$ ($Y_k=Y_\pm$). Condition~\ref{cond:trans} guarantees that the
flow $Y_k$ intersects the switching manifold $\{|h|=\epsilon\}$
transversally in $y_*(t_{*,k})$ for all crossing times $t_{*,k}$, that is,
$h'(y_*(t_{*,k}))\dot Y_k^t(y_*(t_{*,k}))\vert_{t=0}\neq 0$.

The coordinates $y_0\in\Sigma$ and $(t_1,\ldots,t_\mu)$ of an initial
condition in $(\xi,1)\in{\cal S}_1$ must be in a small neighborhood of
$(y_*(0),t_{*,m-\mu+1}-p,\ldots,t_{*,m}-p)$. Therefore, the headpoint
$y(t)$ of $E_c^t(\xi,1)$ follows $Y_k$ at times near $t_{*,k}$ for
$k=1,\ldots,m$ (because $y_*$ follows $Y_k$). Thus, the transversal
intersection of $Y_k$ with the switching manifold near
$t_{*,1},\ldots,t_{*,m}$ implies that the crossings of $y(t)$ near
$t_{*,1},\ldots,t_{*,m}$ depend smoothly on the coordinates
$(y_0,t_1,\ldots,t_\mu)$.  Consequently, also all times $t\in(0,p)$
when $E_d^t(\xi,1)$ changes its value depend smoothly on
$(y_0,t_1,\ldots,t_\mu)$. As the cross-section $\Sigma$ is transversal
to $Y_+$ and $y(t)$ follows $Y_+$ for $t$ near $p$ the return time
also depends smoothly on $(y_0,t_1,\ldots,t_\mu)$ and, hence, the
whole return map $P$ depends smoothly on the coordinates
$y_0\in\Sigma$ and $(t_1,\ldots,t_\mu)$.

\section{Proof of Lemma~\ref{thm:retmapreduce} and
  Theorem~\ref{thm:retmapform}}
\label{sec:app:switch}

\paragraph*{Proof of Lemma~\ref{thm:retmapreduce}}
The strict transversality condition \eqref{eq:stran:qg0} implies that
the periodic orbit also satisfies the weak transversality condition
(as given in Definition~\ref{def:wtran}). In the proof of continuity
(Lemma~\ref{thm:continuity}) we established that the discrete
components $E_d^t(\xi,u)$ of trajectories starting from initial
conditions $(\xi,u)$ near $(y_*(\cdot),u_*(0))$ change their values
always close to times where $u_*$ changes its value. Also the
direction of change must be the same because there is a minimum
distance between subsequent changes (given in \eqref{eq:mindist}). The
colliding symmetric periodic orbit changes its value exactly twice per
period, at $t=0$ and at $t=T$. Thus, the image $(\xi,u)$ of an initial
value $(\xi_0,u_0)$ must have a continuous component $\xi$ switching
exactly once near $s=-\Delta$ from $Y_-$ to $Y_+$.  (The Poincar{\'e}
section $\Sigma$ was taken at $y_*(\Delta)$ and $u_*$ switches from
$-1$ to $1$ at $t=0$.)

Consequently, $\xi$ must have the form \eqref{eq:mapxi} for some $y\in
U(y_*(0))$ and a time $\theta(y)$, which is the traveling time from
$y$ to $\Sigma$ following $Y_+$.
The only open question is if 
$\xi(0)\in\Sigma$ and the time $t_1(\xi)$ are uniquely determined by $y$.  

Let $y\in U(y_*(0))$ be given. The time $\theta(y)$ is
implicitly given by $Y_+^{\theta}y\in\Sigma$, which means
\begin{equation}
  f_2^T[Y_+^\theta  y-y_*(\Delta)]=0\mbox{.}\label{eq:thetadef}
\end{equation}
The linearization of \eqref{eq:thetadef} in $y=y_*(0)$ with respect to
$\theta$ is $f_2^Tf_2$, which is non-zero due to the strict
transversality \eqref{eq:stran:qg0}. Thus, $\theta$ is a locally
unique and well defined function of $y$. Consequently,
$\xi(0)=Y_+^{\theta}y$ is also well defined and smoothly dependent on
$y$. The time $t_1=t_1(\xi)$ is given by $t_1=t(y)-\theta(y)$ where
$t$ is implicitly defined by
\begin{eqnarray}
    \epsilon&=&h(Y_+^{t}y)\mbox{\quad if $h(y)<\epsilon$,}
    \label{eq:t1defp}\\
    \epsilon&=&h(Y_-^{t}y)\mbox{\quad if $h(y)\geq\epsilon$.}\label{eq:t1defm}
\end{eqnarray}
The linearization with respect to $t$ in $y=y_*(0)$ is $h'^Tf_2$ in
case \eqref{eq:t1defp} and $h'^Tf_1$ in case \eqref{eq:t1defm}. Both
linearizations are non-zero due to the  strict transversality condition
\eqref{eq:stran:qg0}.

\paragraph*{Proof of Theorem~\ref{thm:retmapform}}
Let $y\in U(y_*(0))$ be given. Lemma~\ref{thm:retmapreduce} gives a
unique element $\xi$ of the image of $P$ corresponding to $y$ that has
the form \eqref{eq:mapxi}. The image of $y$ under the map $m$ is the
location $y_2$ at the next time $s_2$ when the discrete component
$E_d^{s_2}(\xi,1)$ changes its value from $-1$ to $+1$.  The full
reflection symmetry of the periodic orbit $L$ implies that $m$ is the
second iterate of a map $F:U(y_*(0))\mapsto U(y_*(0))$ which is
defined as $-y_1$ where $y_1$ is the location at the next time $s_1$
when the discrete component $E_d^{s_1}(\xi,1)$ changes its value from
$+1$ to $-1$.

The regular implicit condition \eqref{eq:thetadef} defines $\theta(y)$
and $\xi(0)$. The regular implicit condition \eqref{eq:ftdef} (same as
\eqref{eq:t1defp},\,\eqref{eq:t1defm}) defines the time
$t_1(y)=t(y)-\theta(y)$ locally uniquely.  Thus, the $E_d^s(\xi,1)$
changes its value to $-1$ at time $s_1=\tau+t_1(y)$. The point $y_1$
is the headpoint of the continuous component $E_c^{s_1}(\xi,1)$. It
has the form
\begin{displaymath}
  y_1=Y_+^{s_1}\xi(0)=Y_+^{s_1+\theta(y)}y=Y_+^{\tau+t(y)}y  
\end{displaymath}
because $t_1=t(y)-\theta(y)$ where $t(y)$ is given by the regular
implicit condition \eqref{eq:ftdef}. Thus, $F(y)=-y_1$ has the form
claimed in Theorem~\ref{thm:retmapform}.
\end{appendix}
\end{document}